\documentclass[aop,seceqn,citesort,MSNbibl,dvips]{arximspdf}
\usepackage[mathscr]{eucal}
\usepackage{accents}

\doi{10.1214/09-AOP516}
\volume{38}
\issue{4}
\pubyear{2010}
\firstpage{1532}
\lastpage{1569}

\makeatletter

\newcommand{\iint}{\int\!\!\int}
\newcommand{\auf}{[\![}
\newcommand{\zu}{]\!]}

\newtheorem{theorem}{Theorem}[section]
\newtheorem{lem}{Lemma}[section]

\newtheorem{cor}{Corollary}[section]
\newtheorem{pr}{Proposition}[section]

\newproclaim{df}{Definition}[section]

\newproclaim{ex}{Example}[section]
\newproclaim{rem}{Remark}[section]

\makeatother

\begin{document}
\begin{frontmatter}

\title{Stochastic calculus over symmetric Markov processes without time reversal}
\runtitle{Stochastic calculus without time reversal}

\begin{aug}
\author[A]{\fnms{Kazuhiro} \snm{Kuwae}\corref{}\thanksref{t1}\ead[label=e1]{kuwae@gpo.kumamoto-u.ac.jp}}
\runauthor{K. Kuwae}
\affiliation{Kumamoto University}
\address[A]{Department of Mathematics\\
\quad and Engineering\\
Graduate School of Science\\
\quad and Technology\\
Kumamoto University\\
Kumamoto, 860-8555 \\
Japan\\
\printead{e1}} 
\end{aug}

\thankstext{t1}{Supported in part by a Grant-in-Aid for Scientific
Research (C) No. 19540220 
from Japan Society for the Promotion of Science.}

\received{\smonth{3} \syear{2009}}
\revised{\smonth{12} \syear{2009}}

%
\begin{abstract}
We refine stochastic calculus for symmetric Markov processes without
using time reverse operators. Under some conditions on the jump
functions of locally square integrable martingale additive functionals,
we extend Nakao's divergence-like continuous additive functional of
zero energy and the stochastic integral with respect to it under the
law for quasi-everywhere starting points, which are refinements of the
previous results under the law for almost everywhere starting points.
This refinement of stochastic calculus enables us to establish a
generalized Fukushima decomposition for a certain class of functions
locally in the domain of Dirichlet form and a generalized It\^o
formula.
\end{abstract}

%

\pdfkeywords{31C25, 60J25, 60J45, 60J75, Symmetric Markov process,
Dirichlet form, Revuz measure, martingale additive functionals of finite
energy, continuous additive functional of zero energy,
Nakao's CAF of zero energy, Fukushima decomposition,
semi-martingale, Dirichlet processes, stochastic integral,
Ito integral, Fisk--Stratonovich integral,
time reversal operator, dual predictable projection}

\end{frontmatter}

\section{Introduction}
In this paper, under the framework of general symmetric Markov
processes without using time reverse operators,
we give a refinement of stochastic calculus developed in the previous
joint paper \cite{CFKZStoch}.
More precisely, we establish stochastic integrals both of It\^o-type
and of Fisk--Stratonovich-type by
Dirichlet processes by extending
the Nakao's divergence-like continuous additive functional 
of zero energy to a continuous additive functional locally of zero
energy for
a class of locally square integrable martingale additive functionals.
Throughout this paper, we use the terminology \textit{Dirichlet process}
specifically for
an additive functional decomposed into
the sum of a locally square integrable martingale additive functional
and a continuous additive functional (locally) of zero energy, which is not
necessarily a semi-martingale in general; indeed, the notion of
Dirichlet process in a more
general context was introduced by
F\"ollmer \cite{FolDir}. As in \cite{FOT}, stochastic integrals are
defined to be
additive functionals admitting exceptional sets. So all formulas in
this paper
can be regarded as a decomposition of additive functional, which holds
for all time (or up to the life time) with
probability $1$ starting from quasi-everywhere point.

Hereafter, we use the abbreviation CAF (resp., MAF) for continuous
additive functional (resp.,
martingale additive functional).
For a Dirichlet process given by Fukushima's decomposition,
Nakao \cite{Na} defined stochastic integrals integrated by
his divergence-like CAF of zero energy, which
enables us to construct an It\^o-type stochastic integral by the
Dirichlet process.
He also defined
a Fisk--Stratonovich-type integral for symmetric diffusion
processes with no inside killing in order to establish the
stochastic line integral along 1-forms for symmetric diffusion processes
over smooth manifolds and gave an application of stochastic line integral
to a homogenization problem.

On the other hand, Lyons and Zheng \cite{LZdecomposition} and Lyons
and Zheng \cite{LZdecom}
introduced the notion of Fisk--Stratonovich-type integrals in terms of
the sum of forward and backward martingales, which is described by
time reverse operators in the framework of symmetric conservative
diffusion processes.
They proved that their Fisk--Stratonovich-type integrals
are consistent with Nakao's one under the law ${\mathbb{P}}_m$.

In the joint paper \cite{CFKZStoch}, we extend Nakao's divergence-like
CAF of
zero energy in terms of time reverse operators
and define a stochastic integral integrated by this extended CAF
under some mild conditions, which plays an important role in deducing
the perturbation of general symmetric
Markov processes, that is, the combination of the Feynman--Kac formula
and the Girsanov formula (see \cite{CFKZPert,CFKZGenPert});
however, still described under the law ${\mathbb{P}}_m$ except a
special case.

We extend Nakao's CAFs of zero energy and stochastic integrals with
respect to it for more
general integrand and integrator in terms of the the space locally in
the Dirichlet space and
a subclass of locally square integrable MAF on $\auf 0,\zeta\auf $
(Definition \ref{df:localizationGamma}).
We will
define both the It\^o-type and the Fisk--Stratonovich-type stochastic integrals
integrated by (not necessarily continuous)
Dirichlet processes under the law ${\mathbb{P}}_x$ for quasi-everywhere
starting point $x\in E$,
which are described in terms of a subclass of locally square integrable
MAF on $\auf 0,\zeta\auf $
over general symmetric Markov processes (Definitions \ref{df:FiskStratonovich}
and \ref{df:extenstochinteg}).
Our definitions of Fisk--Stratonovich-type integrals are somewhat
different from what is defined by Meyer \cite{Meyerstochs} and Protter
\cite{ProtStoch} in
the framework of semi-martingales (Remark \ref{rem:FiskStratonovich}).

We further show that our stochastic integrals integrated by
the purely discontinuous part of Dirichlet processes
have a representation
of sum of jumps on Dirichlet processes if the jump function of
integrator is anti-symmetric, which
enables us to see the pathwise behavior of pure jump processes under
the law for quasi-everywhere starting points
(Theorem \ref{thm:Stratonovichlocalrepre}, Corollary \ref
{cor:purejumpapprox2}).

As a corollary, we establish a generalized Fukushima decomposition for
a class of functions
locally in the domain of forms (Theorem \ref{thm:FDlocallyinDir}).
We also
present a generalized It\^o formula in terms of our extended stochastic
integrals by
Dirichlet processes (Theorem \ref{thm:ItofornulaFS}).
Our It\^o formula for Fisk--Stratonovich-type integrals
has an expression different from what is exposed in
Protter \cite{ProtStoch} (Remark \ref{rem:ItofornulaFS}).

Let us briefly outline the organization of this paper. In Section \ref
{sec:preliminary}, we
describe the setting of the paper and give some basic lemmas.
In Section \ref{sec:Nakao}, we formulate the extension of Nakao's CAF
of zero energy and
stochastic integral with respect to it under the law for
quasi everywhere starting points.
In Section \ref{sec:StochasticIntegral},
we define our stochastic integrals by Dirichlet processes and
expose the result as noted
above.

\section{Preliminary facts}\label{sec:preliminary}

Let $\mathbf{M}=\{ \Omega, \mathscr{F}_{\infty},\mathscr{F}_t,
X_t,\theta_t, \zeta, {\mathbb{P}}_x, x\in
E\}$ be an $m$-symmetric right Markov process on a Lusin space $E$,
where $m$ is a $\sigma$-finite measure with full support on $E$. Its
associated Dirichlet space $(\mathcal{E}, \mathcal{F})$ on $L^2(E;m)$
is known to be
quasi-regular (see \cite{MR}). By \cite{CMR}, $(\mathcal{E},
\mathcal{F})$ is
quasi-homeomorphic to a regular Dirichlet space on a locally compact
separable metric space.
Thus using this quasi-homeomorphism, without
loss of generality, we may and do assume that \textbf{M} is an
$m$-symmetric Hunt process on a locally compact metric space $E$
such that its associated Dirichlet space $(\mathcal{E}, \mathcal{F})$
is regular on
$L^2(E;m)$ and that $m$ is a positive Radon measure with full
topological support on $E$. But we implicitly use the quasi-left
continuity up to
$\infty$, which is not the usual property of right Markov processes.
So the strict
quasi-regularity of $(\mathcal{E}, \mathcal{F})$ is essentially assumed.
However, if we restrict ourselves to state the result that holds up to
the life time with probability $1$ for quasi-everywhere starting point,
then the framework of quasi-regular Dirichlet forms
is enough.

Without loss of generality, we can take $\Omega$ to be the canonical
path space $D([0,\infty[\ \to E_{\partial})$ of right-continuous,
left-limited (\textit{rcll}, for short) functions from $[0,\infty[$
to $E_{\partial}$. For any $\omega\in\Omega$, we set
$X_t(\omega):=\omega(t)$. Let $\zeta(\omega):=\inf\{t\geq0 \mid
X_t(\omega)=\partial\}$ be the life time of \textbf{M}. As usual,
$\mathscr{F}_{\infty}$ and $\mathscr{F}_t$ are the minimal completed
$\sigma$-algebras
obtained from $\mathscr{F}_{\infty}^0:=\sigma\{X_s \mid  0\leq
s<\infty\}$
and $\mathscr{F}_t^0:=\sigma\{X_s \mid  0\leq s\leq t\}$, respectively,
under ${\mathbb{P}}_x$.
We set
$X_t(\omega):=\partial$ for $t\geq\zeta(\omega)$ and use $\theta_t$
to denote the shift operator defined by
$\theta_t(\omega)(s):=\omega(t+s)$, $t,s\geq0$.
For each $s\geq0$, the shift operator $\theta_s$ is defined by
$\theta_s\omega(t):=\omega(t+s)$ for $t\in[0,\infty[$.
For a Borel subset $B$ of $E$,
$\sigma_B:=\inf\{t>0\mid X_t\in B\}$ (the first hitting time to $B$) and
$\tau_B:=\inf\{t>0\mid X_t\notin
B\}$ (the \textit{first exit time} of $B$) are $(\mathscr
{F}_t)$-stopping times. If
$B$ is closed, then $\tau_B$ is an $(\mathscr{F}_{t+}^0)$-stopping time.
Also, $\zeta$ is an $(\mathscr{F}_t^0)$-stopping time because $\{
\zeta\leq
t\}=\{X_t=\partial\}\in\mathscr{F}_t^0$, $t\geq0$.

The transition semigroup of \textbf{M}, $\{P_t, t\ge0\}$, is defined by
\[
P_tf(x):={\mathbb{E}}_x[f(X_t)]={\mathbb{E}}_x[f(X_t)\dvtx t< \zeta
],\qquad  t\ge0.
\]
Each $P_t$ may be viewed as an operator on $L^2(E;m)$; collectively
these operators form a strongly continuous semigroup of
self-adjoint contractions. The Dirichlet form associated with \textbf{M}
is the
bilinear form
\[
{\mathcal{E}}(u, v):=\lim_{t\downarrow0}\frac1t(u-P_tu, v)_m
\]
defined on the space
\[
{\mathcal{F}}:= \Bigl\{u\in L^2(E;m)  \bigm|
 \sup_{t>0}  t^{-1}(u-P_tu, u)_m<\infty\Bigr\}.
\]
Here we use the notation $(f,g)_m:=\int_E f(x)g(x) m(dx)$ for $f,g\in
L^2(E;m)$.

An increasing sequence $\{F_n\}$ of closed sets
is called an \textit{$\mathcal{E}$-nest} if $\bigcup_{n=1}^\infty
\mathcal{F}_{F_n}$ is
$\mathcal{E}_1^{1/2}$-dense in $\mathcal{F}$, where $\mathcal
{F}_{F_n}:=\{u\in\mathcal{F}\mid  u=0\mbox{ }
m\mbox{-a.e.
on } E\setminus F_n\}$ and a family $\{F_n\}$ of closed sets is an
$\mathcal{E}
$-nest if and only
if it is a nest, that is, ${\mathbb{P}}_x(\lim_{n\to\infty}\tau
_{F_n}=\zeta
)=1$ q.e. $x\in E$.
A function $u$ on $E$ is said to be \textit{$\mathcal{E}$-quasi-continuous}
if there exists an $\mathcal{E}$-nest $\{F_n\}$
of closed sets such that $u$ is continuous on each $F_n$.
A subset $N$ of $E$ is called \textit{$\mathcal{E}$-polar or} ($\mathcal{E}
$-)\textit{exceptional} if there exists an $\mathcal{E}$-nest
$\{F_n\}$ such that $N\subset\bigcap_{n=1}^{\infty}(E\setminus
F_n)$; equivalently
there is a Borel set $\tilde{N}$ containing $N$ such that ${\mathbb{P}}
_m(\sigma_{\tilde{N}}<\infty)=0$.
A statement $S(x)$ is said to hold for quasi-everywhere $x\in E$
(q.e. $x\in E$ in short) if
there exists an exceptional set $N$ such that
$\{x\in E\mid S(x)$ does not hold$\}\subset N$.

An increasing sequence $\{F_n\}$ of closed sets
is called a \textit{strict $\mathcal{E}$-nest} if
\[
\lim_{n\to\infty}\operatorname{Cap}_{1,G_1\varphi}(E\setminus F_n)=0,
\]
where
$\operatorname{Cap}_{1,G_1\varphi}$ is the weighted capacity defined
in Chapter V, Definition 2.1 of
\cite{MR} and a family $\{F_n\}$ of closed sets is a strict $\mathcal{E}$-nest
if and only
if it is a strict nest, that is, ${\mathbb{P}}_x(\lim_{n\to\infty
}\sigma
_{E\setminus F_n}=\infty)=1$
$m$-a.e. $x\in E$ in view of Chapter V, (2.5) in~\cite{MR},
equivalently it holds q.e. $x\in E$
by Chapter V, Proposition 2.28(i) and Remark 2.8 in \cite{MR}.
A function $u$ on $E_{\partial}$ is said
to be \textit{strictly $\mathcal{E}$-quasi-continuous}
if there exists a strict $\mathcal{E}$-nest $\{F_n\}$
of closed sets such that $u$ is continuous on each $F_n\cup\{\partial
\}$.

An increasing sequence $\{G_n\}$ of (q.e.) finely open Borel sets is
called a \textit{nest} (resp., \textit{strict nest}) if
${\mathbb{P}}_x(\lim_{n\to\infty}\tau_{G_n}=\zeta)=1$ for
q.e. $x\in E$ [resp.,
${\mathbb{P}}_x(\lim_{n\to\infty}\sigma_{E\setminus G_n}=\infty
)=1$ for
q.e. $x\in E$].
(The definition of q.e. finely open sets can be found in \cite{FOT}.)
In \cite{CFKZStoch}, we show that under the quasi-left-continuity up
to infinity of \textbf{M},
for an increasing sequence $\{G_n\}$ of (q.e.) finely open Borel sets,
$\{G_n\}$ is a nest
if and only if it is a strict nest. Denote by $\Theta$ the family of
(strict) nests $\{G_n\}$ of (q.e.) finely open Borel sets.
Note that
for an $\mathcal{E}$-nest $\{F_n\}$ of closed sets, $\{G_k\}\in\Theta
$ by setting
$G_k:=F_k^{f\mbox{-}\mathrm{int}}$, $k\in\mathbb{N}$, where
$F_k^{f\mbox{-}\mathrm{int}}$ means the fine interior of $F_k$.

Let $\mathcal{F}_e$ be the family of $m$-measurable functions $u$ on $E$
such that $\vert u\vert<\infty$ $m$-a.e. and there exists an
${\mathcal{E}}$-Cauchy
sequence $\{u_n\}$ of $\mathcal{F}$ such that $\lim _{n\to
\infty}u_n
=u m$-a.e.
We call $\{u_n\}$ as above an approximating sequence for $u\in\mathcal{F}_e$.
For any $u,v\in\mathcal{F}_e$ and their approximating sequences
$\{u_n\}$, $\{v_n\}$
the limit ${\mathcal{E}}(u,v)=\lim _{n\to\infty}{\mathcal
{E}}(u_n,v_n)$
exists and does not
depend on the choices of the approximating sequences for $u$, $v$.
It is known that ${\mathcal{E}}^{1/2}$ on $\mathcal{F}_e$ is a
semi-norm and $\mathcal{F}=\mathcal{F}
_e\cap L^2(E;m)$.
We call $(\mathcal{E},\mathcal{F}_e)$ \textit{the extended Dirichlet
space} of $(\mathcal{E},\mathcal{F})$.
Let $L^0(E;m)$ be the family of $m$-measurable functions on $E$.
We further let
\begin{eqnarray*}
\dot{\mathcal{F}}_{\mathrm{loc}}&:=&\bigl\{u\in L^0(E;m)\mid
\mbox{ there exist }\{G_n\}\in\Theta\mbox{ and }\\
&&\hspace*{5.1pt}{u}_n\in\mathcal
{F}\mbox{
such that } u=u_n\mbox{ }
m\mbox{-a.e. on }G_n\mbox{ for each }n\in\mathbb{N}\bigr\}.
\end{eqnarray*}
$\dot{\mathcal{F}}_{\mathrm{loc}}$ is called the \textit{space of functions
locally in
$\mathcal{F}$ in the broad sense}.
It is shown in \cite{Kwfunc} that $\mathcal{F}\subset\mathcal
{F}_e\subset\dot{\mathcal{F}
}_{\mathrm{loc}}$ and every
$u\in\dot{\mathcal{F}}_{\mathrm{loc}}$ admits an $\mathcal
{E}$-quasi-continuous $m$-version
$\tilde{u}$. More strongly
every $u\in\mathcal{F}$ admits a strictly $\mathcal
{E}$-quasi-continuous $m$-version
$\tilde{u}$ on $E_{\partial}$ with $\tilde{u}(\partial)=0$. For
$u\in\dot{\mathcal{F}}_{\mathrm{loc}}$, we always assume that $\tilde{u}$ is extended
to be a
real-valued function $\tilde{u}$ on $E_{\partial}$ if not otherwise
specified, where we do not necessarily assume $\tilde{u}(\partial)=0$.
However, we can reduce to this case by setting
$\tilde{u}-\tilde{u}(\partial)$
on $E_{\partial}$, which is in $\dot{\mathcal{F}}_{\mathrm{loc}}$ as a function
defined on $E$.

We need the following lemma:
\begin{lem}\label{lem:strictqc}
Every $u\in\mathcal{F}_e$ admits a strictly $\mathcal
{E}$-quasi-continuous $m$-version
$\tilde{u}$ on $E_{\partial}$ with $\tilde{u}(\partial)=0$.
\end{lem}
\begin{pf}
Take $u\in\mathcal{F}_e$. Then there exists an $m$-a.e. strictly positive
bounded function $g\in L^1(E;m)$
such that $u\in(\mathcal{F}^g)_e$, where $(\mathcal{E}^g,\mathcal
{F}^g)$ is the Dirichlet form
on $L^2(E;m)$ defined
by $\mathcal{F}^g:=\mathcal{F}\cap L^2(E;gm)$, $\mathcal
{E}^g(v,w):=\mathcal{E}(v,w)+(v,w)_{gm}$,
$v,w\in\mathcal{F}^g$, and
$(\mathcal{F}^g)_e$ is its extended Dirichlet space.
Then there exist an increasing sequence $\{F_n\}$ of closed sets and a
function $\tilde{u}$
on $E_{\partial}$
such that
\[
\lim_{n\to\infty}\operatorname{Cap}_{(0)}^g(E\setminus F_n)=0,
\]
$\tilde{u}=u$ $m$-a.e. on $E$
and $\tilde{u}$ is continuous on each $F_n\cup\{\partial\}$ with
$\tilde{u}(\partial)=0$,
where $\operatorname{Cap}_{(0)}^g$ is the $0$-order capacity with
respect to $(\mathcal{E}^g,(\mathcal{F}^g)_e)$.
It suffices to prove that $\{F_n\}$ is a strict $\mathcal{E}$-nest
with respect
to $(\mathcal{E},\mathcal{F})$.
For this, we need that for any open set $U$,
\[
H_U^g1(x):={\mathbb{E}}_x \bigl[e^{-\int_0^{\sigma_U}g(X_s)\,ds} \bigr]
\]
satisfies $H_U^g1\in(\mathcal{F}^g)_e$ and
%
\begin{equation}\label{eq:equibilium0order}
\operatorname{Cap}_{(0)}^g(U)=\mathcal{E}^g(H_U^g1,H_U^g1).
\end{equation}
This can be similarly proved along the same way as in
Section 4.4 in \cite{FOT}. We will omit the details.

From (\ref{eq:equibilium0order}), we have
\[
{\mathbb{P}}_x \biggl(\lim_{n\to\infty}
\int_0^{\sigma_{E\setminus F_n}}g(X_s)\,ds=\infty\biggr)=1,\qquad  m\mbox
{-a.e. }x\in E,
\]
and hence $\{F_n\}$ is a strict nest, because of the boundedness of $g$.
\end{pf}

As a rule we take $u$ to be
represented by its (strictly) $\mathcal{E}$-quasi-continuous
$m$-version (when
such exists), and
drop the tilde from the notation.

Let $\accentset{\circ}{\mathcal{M}}$ and $\mathcal{N}_c$ denote,
respectively, the
space of martingale additive functionals
of finite energy and the space of continuous additive
functionals
of zero energy.
More precisely, we set
\begin{eqnarray*}
\mathcal{M}:\!&=&\{M\mid M\mbox{ is a finite rcll AF},  {\mathbb
{E}}_x[M_t^2]<\infty,
 {\mathbb{E}}_x[M_t]=0 \\
&&\hspace*{111.8pt}
\mbox{ for q.e. } x\in E\mbox{ and all }t\geq0\}.
\end{eqnarray*}
For an AF $M$, if the limit
%
\begin{equation}
\mathbf{e}(M):=\lim_{t\downarrow0}\frac1{2t}{\mathbb{E}}_m[M_t^2]
\end{equation}
exists, we call it \textit{energy of} $M$. When $M\in\mathcal{M}$,
$t\mapsto
\frac1{2t}{\mathbb{E}}_m[M_t^2]$ is
increasing and the limit may diverge in general. Then we define
\begin{eqnarray*}
\accentset{\circ}{\mathcal{M}}&:=& \{M\in{\mathcal{M}}  \mid
\mathbf{e}(M)
<\infty\}, \\
\mathcal{N}_c&:=& \{N  \mid  N\mbox{ is a finite CAF},
{\mathbb{E}}_x[|N_t|]<\infty\mbox{ q.e. } x\in E \\
&&\hspace*{84.6pt}
\mbox{ for
each }t>0\mbox{,
 and }\mathbf{e}(N)=0 \}.
\end{eqnarray*}
For $M,N\in\accentset{\circ}{\mathcal{M}}$,
we set
\[
\mathbf{e}(M,N):=\lim_{t\downarrow0}\frac1{2t}{\mathbb{E}}_m[M_tN_t]
\]
and call it \textit{mutual energy} of $M,N$. It is well known that
$(\accentset{\circ}{\mathcal{M}},\mathbf{e})$ is a real Hilbert
space with
inner product $\mathbf{e}$.

For $u\in\mathcal{F}_e$, the following Fukushima
decomposition holds:
%
\begin{equation}\label{eqn:fuku}
u(X_t)-u(X_0)=M^u_t+N^u _t
\end{equation}
for all $t\in[0,\infty[$ ${\mathbb{P}}_x$-a.s. for q.e. $x\in E$,
where $M^u \in
\accentset{\circ}{\mathcal{M}}$ and $N^u\in\mathcal{N}_c$.

A positive continuous additive functional (PCAF) of \textbf{M} (call it $A$)
determines
a measure $\nu=\nu_A$ on the Borel subsets of $E$ via the formula
%
\begin{equation}\label{eq:revuz}
\nu(f)=\uparrow\lim_{t\to0}\frac1t {\mathbb{E}}_m \biggl[\int
_0^tf(X_s) \,dA_s \biggr],
\end{equation}
in which $f\dvtx E\to[0,\infty]$ is Borel measurable. The measure $\nu$
is necessarily \textit{smooth} (denote by $\nu\in S$), in the sense
that $\nu$ charges no
exceptional set of $E$, and there is an $\mathcal{E}$-nest $\{F_n\}$
of closed
subsets of $E$ such that $\nu(F_n)<\infty$ for each
$n\in\mathbb{N}$.
Conversely, given a smooth measure $\nu$,
there is a unique PCAF $A^\nu$ such that (\ref{eq:revuz}) holds
with $A=A^\nu$. In the sequel we refer to this bijection between
smooth measures and PCAFs as the \textit{Revuz correspondence}, and to
$\nu$ as the Revuz measure of~$A^\nu$.
\begin{lem}\label{lem:smoothfinelyopen}
$\mu\in S$ if and only if $\mu$ charges no exceptional set and there
exists $\{G_n\}\in\Theta$ such that
$\mu(G_n)<\infty$ for each $n\in\mathbb{N}$.
\end{lem}
\begin{pf}
The ``only if'' part is trivial by setting $G_n:=F_n^{f\mbox{-}\mathrm{int}}$.
We only prove the ``if'' part.
Take an $f\in L^2(E;m)$ with $0<f\leq1$ on $E$ and set
$R_1^{G_n}f(x):={\mathbb{E}}_x[\int_0^{\tau_{G_n}}e^{-s}f(X_s)\,ds]$. Then
$R_1^{G_n}f(x)>0$ on $G_n$ and $R_1^{G_n}f$ is $\mathcal{E}$-quasi-continuous
for each $n\in\mathbb{N}$. Take a common $\mathcal{E}$-nest $\{A_k\}
$ such that
all $R_1^{G_n}f$, $n\geq1$ are continuous on each $A_k$.
We set $F_n:=\{x\in A_n\mid R_1^{G_n}f(x)\geq1/n\}$. Then $\{F_n\}$ is
an $\mathcal{E}$-nest
by use of Lemma 3.3 in \cite{Kwfunc}, where we observe $B_n:=\{
R^{G_n}_1f>1/n\}$ is increasing
and $E\setminus\bigcup_{n=1}^{\infty}B_n$ is exceptional.
For each $n\in\mathbb{N}$, we have $(E\setminus G_n)^r\subset
E\setminus F_n$, where
$(E\setminus G_n)^r=\{x\in E\mid R_1^{G_n}f(x)=0\}$ is the set of
regular points for
$E\setminus G_n$. Since $(E\setminus G_n)\setminus(E\setminus G_n)^r$ is
exceptional, we obtain $\mu(F_n)\leq\mu(G_n)<\infty$ for each $n\in
\mathbb{N}$.
\end{pf}

A (positive) Radon measure $\mu$ on $E$ is said to be a measure
\textit{of finite energy integral} if there exists $C>0$
depending on $\mu$ such that
\[
\int_E|u(x)|\mu(dx)\leq C\sqrt{\mathcal{E}_1(u,u)}\qquad\mbox{for all
}u\in\mathcal{F}
\cap C_0(E).
\]
Let $S_0$ be the family of measures of finite energy integrals. For
$\mu\in S_0$ and $\alpha>0$,
there exists a unique element $U_{\alpha}\mu\in\mathcal{F}$ such that
\[
\mathcal{E}_{\alpha}(U_{\alpha}\mu,v)=\int_Ev(x)\mu(dx)\qquad
\mbox{for }v\in
\mathcal{F}\cap C_0(E).
\]
It is known that every $\mu\in S_0$ is a smooth measure.
If we set $S_{00}:=\{\mu\in S_0\mid\mu(E)<\infty, U_1\mu\in
L^{\infty}(E;m)\}$, then
$N$ is exceptional if and only if $\nu(N)=0$ for all $\nu\in S_{00}$.

For any $\mu\in S$, $\nu\in S_{00}$, a (q.e.) finely open Borel set
$G$ and $t>0$,
we have the following formula:
%
\begin{equation}\label{eq:fineversionoflem5.1.9}
{\mathbb{E}}_{\nu}[A^{\mu}_{t\wedge\sigma_{E\setminus G}}]\leq
(1+t)\|
U_1\nu\|_{\infty}
\mu(G),
\end{equation}
which can be similarly proved as in the proof of Lemma 5.1.9 in \cite{FOT}
with the help of Lemma 5.1.10(ii) in \cite{FOT}.

Take $M,N\in\accentset{\circ}{\mathcal{M}}$ and
denote by $\langle M,N\rangle$ its quadratic covariational process,
which is a
CAF of bounded variation, and let $\mu_{\langle M,N\rangle}$ be its
Revuz measure.
In view of Theorem 2.2 in \cite{Na}, for $M,N\in\accentset{\circ}{\mathcal{M}}$,
$\mathbf{e}(M,N)=0$ implies
that $\langle M,N\rangle\equiv0$ on $[0,\infty[$ ${\mathbb{P}}_x$-a.s. for
q.e. $x\in E$.
For $M\in\accentset{\circ}{\mathcal{M}}$ and $f\in L^2(E;\mu
_{\langle
M\rangle})$, there exists a unique $f*M\in\accentset{\circ}{\mathcal{M}}$
such that
\[
\mathbf{e}(f*M,N)=\frac12 \int_Ef(x)\mu_{\langle M,N\rangle
}(dx)\qquad\mbox{for }N\in\accentset{\circ}{\mathcal{M}}.
\]
Moreover, we have the following.
\begin{lem}\label{lem:Riemannsumapproxmar}
Let $M\in\accentset{\circ}{\mathcal{M}}$ and $f\in L^2(E;\mu
_{\langle
M\rangle})$.
If $f$ is a strictly $\mathcal{E}$-quasi-continuous function, then
$f*M$ admits a
Riemann sum approximation: for each $t>0$
\[
(f*M)_t=\lim_{n\to\infty}\sum_{\ell=0}^{n-1}f(X_{\ell
t/n})\bigl(M_{(\ell+1)t/n}-M_{\ell t/n}\bigr)
\]
holds ${\mathbb{P}}_x$-a.s. for q.e. $x\in E$, where the convergence
of the
right-hand side is
in ${\mathbb{P}}_x$-probability for q.e. $x\in E$.
\end{lem}
\begin{pf}
This is well known for experts and shown for the case $f\in C_0(E)$ in
Lemma 5.6.2 of \cite{FOT}.
We shall show it for the reader's convenience. By assumption, we have
for $\nu\in S_{00}$
\[
{\mathbb{E}}_{\nu} \biggl[\int_0^tf^2(X_s)\,d\langle M\rangle_s \biggr]\leq
(1+t)\|U_1\nu
\|_{\infty}\int_Ef^2(x)\mu_{\langle M\rangle}(dx)<\infty.
\]
In particular,
\[
{\mathbb{E}}_x \biggl[\int_0^tf^2(X_s)\,d\langle M\rangle_s \biggr]<\infty \qquad\mbox
{for q.e. }x\in E.
\]
Then by Theorem A.3.19 in \cite{FOT}, for $x\in E\setminus N$ with an
adequate properly exceptional set $N$,
we can define the stochastic integral $f\bullet M:=\int
_0^tf(X_{s-})\,dM_s$ under~${\mathbb{P}}_x$,
which is characterized by
\[
{\mathbb{E}}_x[(f\bullet M)_t^2]={\mathbb{E}}_x \biggl[\int
_0^tf^2(X_{s-})\,d\langle M\rangle
_s \biggr]
={\mathbb{E}}_x \biggl[\int_0^tf^2(X_s)\,d\langle M\rangle_s \biggr].
\]
From this, we can get $f\bullet M \in\accentset{\circ}{\mathcal{M}}$ and
$\mu_{\langle f\bullet M,N\rangle}=
f\mu_{\langle M,N\rangle}$ for $N\in\accentset{\circ}{\mathcal{M}}$, hence
we have $f*M=f\bullet M$. On the other hand, since $t\mapsto f(X_{t-})$ is
left-continuous ${\mathbb{P}}_x$-a.s. for q.e. $x\in E$,
$f\bullet M$ admits the Riemann-sum approximation in view of Chapter I,
Proposition 4.44 in
\cite{JS}. Therefore, we obtain the result.
\end{pf}
\begin{rem}
From Lemma \ref{lem:Riemannsumapproxmar},
we may write $(f*M)_t=\int_0^t f(X_{s-})\,dM_s$
if $f$ is strict $\mathcal{E}$-quasi-continuous on $E_{\partial}$.
\end{rem}

Let $(N(x, dy), H_t)$ be a L\'evy system for \textbf{M};
that is, $N(x, dy)$ is a kernel
on $(E_{\partial}, {\mathscr{B}}(E_{\partial}))$
and $H_t$ is a PCAF
with bounded $1$-potential such that for any nonnegative
Borel function $\phi$ on $E_{\partial}\times E_{\partial}$
vanishing on the diagonal
and any $x\in E_{\partial}$,
\[
{\mathbb{E}}_x \biggl[\sum_{s\le t} \phi(X_{s-}, X_s) \biggr]
={\mathbb{E}}_x \biggl[\int^t_0\int_{E_{\partial}} \phi(X_s, y)N(X_s,
dy)\,dH_s \biggr].
\]
To simplify notation, we will write
\[
N \phi(x):= \int_{E_{\partial}} \phi(x, y) N(x, dy).
\]
Let $\mu_H$ be
the Revuz measure of the PCAF $H$. Then the jump measure
$J$ and the killing measure $\kappa$ of \textbf{M} are given by
\[
J(dx\,dy)=\tfrac12N(x, dy)\mu_H (dx)\quad\mbox{and}\quad
\kappa(dx)=N(x, \{\partial\})\mu_H(dx).
\]
These measures feature in the Beurling--Deny decomposition of $\mathcal{E}$:
for $f, g\in{\mathcal{F}}_e$,
\begin{eqnarray*}
{\mathcal{E}}(f, g) &=& {\mathcal{E}}^{c}(f, g)+
\int_{E\times E}\bigl(f (x)- f(y)\bigr)
\bigl( g(x)- g(y)\bigr)J(dx\,dy)\\
&&{} +\int_E f(x)
g(x)\kappa(dx),
\end{eqnarray*}
where ${\mathcal{E}}^{c}$ is the strongly local part of ${\mathcal{E}}$.

For $u\in\mathcal{F}_e$, the martingale part $M_t^u$ in (\ref
{eqn:fuku}) can
be decomposed as
\[
M_t^u=M_t^{u,c}+M_t^{u,j} +M_t^{u,\kappa}  \qquad\mbox{for every }
t\in[0,\infty[,
\]
${\mathbb{P}}_x$-a.s. for q.e. $x\in E$, where $M_t^{u,c}$ is the continuous
part of martingale $M^u$, and
\begin{eqnarray*}
M_t^{u,j} & =& \lim_{{\varepsilon}\downarrow0} \biggl\{ \sum_{0<s\leq t}
\bigl( u
(X_s) - u (X_{s-})\bigr)\mathbf{1}_{\{| u (X_s) - u (X_{s-})|>{\varepsilon
}\}}
\mathbf{1}_{\{s<\zeta\}} \\
&&\hspace*{18.8pt}{}  -\int_0^t \biggl( \int_{\{y\in E\mid  |u (y)-
u(X_s)|>{\varepsilon}\}}
\bigl( u (y)- u (X_s)\bigr)
N(X_s,dy) \biggr)\,dH_s \biggr\}, \\
M_t^{u, \kappa}& =& \int_0^t u(X_s) N(X_s,  \{\partial\})\,dH_s-
u(X_{\zeta-})\mathbf{1}_{\{t\geq\zeta\}}
\end{eqnarray*}
are the jump and killing parts of $M^u$ in $\accentset{\circ}{\mathcal{M}}$,
respectively. The limit in the expression for $M^{u, j}$ is
in the sense of convergence in $\accentset{\circ}{\mathcal{M}}$
and of convergence in probability
under ${\mathbb{P}}_x$ for q.e. $x\in E$ for each fixed $t>0$. (See Theorem
A.3.9 and page 341 in \cite{FOT}.)

If we let
\begin{eqnarray*}
\accentset{\circ}{\mathcal{M}}^c&:=&\{M\in\accentset{\circ}{\mathcal{M}}
\mid
M\mbox{ is a continuous MAF}\}, \\
\accentset{\circ}{\mathcal{M}}^d&:=&(\accentset{\circ}{\mathcal{M}}^c)^{\bot}=
\{M\in\accentset{\circ}{\mathcal{M}}\mid\mathbf{e}(M,N)=0\mbox
{ for
}N\in
\accentset{\circ}{\mathcal{M}}^c\},
\end{eqnarray*}
then every $M$ has an orthogonal decomposition
\[
M=M^c+M^d
\]
in the Hilbert space $(\accentset{\circ}{\mathcal{M}},\mathbf
{e})$. $M^c\in
\accentset{\circ}{\mathcal{M}}^c$ (resp., $M^d\in\accentset{\circ}{\mathcal{M}}^d$)
is nothing but the \textit{continuous part} (resp., \textit{purely
discontinuous part}) of $M$
discussed in \cite{FOT}. Moreover, set
\[
\accentset{\circ}{\mathcal{M}}^j:=\{M\in\accentset{\circ}{\mathcal{M}}^d \mid
\mathbf{e}(M,M^{u,\kappa})=0\mbox{ for }u\in\mathcal{F}_e\},\qquad
\accentset{\circ}{\mathcal{M}}^{\kappa}:=\accentset{\circ}{\mathcal{M}}^d\cap
(\accentset{\circ}{\mathcal{M}}^j)^{\bot}.
\]
Then $\accentset{\circ}{\mathcal{M}}^j$ is a closed subspace of
$\accentset{\circ}{\mathcal{M}}$, hence
$M^d$ has a unique orthogonal decomposition in $(\accentset{\circ}{\mathcal{M}}
,\mathbf{e})$ as
\[
M^d=M^j+M^{\kappa},
\]
where $M^j\in\accentset{\circ}{\mathcal{M}}^j$ and $M^{\kappa}\in
\accentset{\circ}{\mathcal{M}}^{\kappa}$.
For simplicity of notation,
we will use the convention $\Delta F_s:=F_s-F_{s-}$ for any rcll
$(\mathscr{F}_t)$-adapted
process $F$.
The square bracket $[M,N]$ for $M,N\in\accentset{\circ}{\mathcal{M}}$ is
defined by
\[
[M,N]_t:=\langle M^c,N^c\rangle_t+\sum_{0<s\leq t}\Delta M_s\Delta N_s.
\]
Then $\langle M,N\rangle$ is the dual predictable projection of
$[M,N]$ (see (A.3.7) in \cite{FOT}).
We further set for each $i=\varnothing,c,d,j,\kappa$
\begin{eqnarray*}
\accentset{\circ}{\mathcal{M}}^i_{f\mbox{-}\mathrm{loc}}:\!&=&\bigl\{M\mid
\mbox{there exist }\{G_n\}\in\Theta\mbox{ and }\bigl\{M^{(n)}\bigr\}\subset
\accentset{\circ}{\mathcal{M}}^i\mbox{ such that }
\\
&&\hspace*{5.1pt}M_t=M_t^{(n)}\mbox{ for all }
t<\sigma_{E\setminus G_n}\mbox{ and }n\in\mathbb{N}, {\mathbb
{P}}_x\mbox
{-a.s. for q.e. }x\in E\bigr\}, \\
{\mathcal{N}}_{c,f\mbox{-}\mathrm{loc}}:\!&=&\bigl\{N\mid\mbox{there exist }\{
G_n\}\in
\Theta\mbox{ and }\bigl\{N^{(n)}\bigr\}\subset\mathcal{N}_c\mbox{ such that }
\\
&&\hspace*{5.1pt}N_t=N_t^{(n)}\mbox{ for all }
t<\sigma_{E\setminus G_n}\mbox{ and }n\in\mathbb{N}, {\mathbb
{P}}_x\mbox
{-a.s. for q.e. }x\in E\bigr\}.
\end{eqnarray*}
Similarly, we can define
$\accentset{\circ}{\mathcal{M}}^i_{\mathrm{loc}}$ and
${\mathcal{N}}_{c,\mathrm{loc}}$ as subclasses of local AFs (or AFs
on $\auf 0,\zeta\auf $)
in terms of first exit times $\tau_{G_n}$ (see \cite{FOT,CFKZStoch}
for the notion of local AF).
Here $i=\varnothing$ means $\accentset{\circ}{\mathcal{M}}^{\varnothing
}:=\accentset{\circ}{\mathcal{M}}$ and
write $\accentset{\circ}{\mathcal{M}}_{f\mbox{-}\mathrm{loc}}$
(resp., $\accentset{\circ}{\mathcal{M}}_{\mathrm{loc}}$)
instead of
$\accentset{\circ}{\mathcal{M}}^{\varnothing}_{f\mbox{-}\mathrm{loc}}$
(resp., $\accentset{\circ}{\mathcal{M}}^{\varnothing
}_{\mathrm{loc}}$).
Every PCAF is an element of $ {\mathcal{N}}_{c,\mathrm{loc}}$.
Our $\accentset{\circ}{\mathcal{M}}_{f\mbox{-}\mathrm{loc}}$
(resp., $ {\mathcal{N}}_{c,f\mbox{-}\mathrm{loc}}$) is
slightly narrower than
$\accentset{\circ}{\mathcal{M}}_{\mathrm{loc}}$ (resp., $
{\mathcal
{N}}_{c,\mathrm{loc}}$)
treated in \cite{FOT} (in \cite{CFKZStoch} we use the same symbol
$\accentset{\circ}{\mathcal{M}}_{f\mbox{-}\mathrm{loc}}$ (resp., $
{\mathcal{N}}_{c,f\mbox{-}\mathrm{loc}}$)
to denote $\accentset{\circ}{\mathcal{M}}_{\mathrm{loc}}$ (resp., $
{\mathcal{N}}_{c,\mathrm{loc}}$)).
However, Fukushima's decomposition (\ref{eqn:fuku}) for $u\in\dot
{\mathcal{F}}_{\mathrm{loc}}$
with $J=\kappa=0$
can be characterized by
$\accentset{\circ}{\mathcal{M}}_{f\mbox{-}\mathrm{loc}}$ and $
{\mathcal
{N}}_{c,f\mbox{-}\mathrm{loc}}$.
Before seeing this, we need the following lemma:
\begin{lem}\label{lem:MAFcontpartvanish}
Let $G$ be a (q.e.) finely open Borel set.
\begin{enumerate}[(1)]
\item[(1)]
If $u\in\mathcal{F}$ satisfies $u=0$
q.e. on $G$, then $\mu_{\langle M^{u,c}\rangle}(G)=0$ and
\[
M_t^{u,c}=0\qquad\mbox{for any }t\leq\sigma_{E\setminus G}\mbox{ }{\mathbb
{P}}_x\mbox
{-a.s. for q.e. }x\in E.
\]
\item[(2)]
For $M\in\accentset{\circ}{\mathcal{M}}$, $\mu_{\langle M\rangle
}(G)=0$ implies
$M_t=0$ for any $t<\sigma_{E\setminus G}$ ${\mathbb{P}}_x$-a.s. for
q.e. $x\in E$.
\end{enumerate}
\end{lem}
\begin{pf}
The proof of (1)
is quite similar to the proof of Lemma 5.3.1 in~\cite{FOT}.
Note that
$t<\sigma_{E\setminus G}\leq\widehat{\sigma}_{E\setminus G}$ implies
$X_s,X_{s-}\in G\cup\{\partial\}$ for all $s\in\ ]0,t]$, which means
$u(X_s)-u(X_{s-})=0$ for all $s\in\ ]0,t]$, because of
$u(\partial)=0$.
Here $\widehat{\sigma}_{E\setminus G}:=\inf\{t>0\mid X_{t-}\in
E\setminus
G\}$ (see (A.2.6) and
Theorem A.2.3 in \cite{FOT}).
Next we prove (2).
Suppose $\mu_{\langle M\rangle}(G)=0$ for $M\in\accentset{\circ}{\mathcal{M}}$.
Note that $\int_0^t\mathbf{1}_{\{\partial\}}(X_s)\,d\langle M\rangle_s=0$
${\mathbb{P}}_x$-a.s. for q.e. $x\in E$. Combining this and
Theorem 5.1.3(i) in
\cite{FOT},
we have ${\mathbb{E}}_m[\int_0^t\mathbf{1}_{G\cup\{\partial\}
}(X_s)\,d\langle
M\rangle_s]=0$ for each $t>0$, hence
$\langle M\rangle_t=0$ for all $t<\sigma_{E\setminus G}$ ${\mathbb
{P}}_m$-a.e.
Then by Lemma 5.1.10(iii) in \cite{FOT}, we obtain the result.
\end{pf}
\begin{rem}
Our method of the proof of Lemma \ref{lem:MAFcontpartvanish}(1)
does not work to show the same assertion in the case that $u$ is only
constant q.e. on $G$.
\end{rem}

From this lemma, we can construct $M^{u,c}\in
\accentset{\circ}{\mathcal{M}}^c_{f\mbox{-}\mathrm{loc}}$ for
$u\in\dot{\mathcal{F}}_{\mathrm{loc}}$.
Under $J=\kappa=0$, for $u\in\dot{\mathcal{F}}_{\mathrm{loc}}$, (\ref{eqn:fuku})
holds for all
$t\in[0,\zeta[$ ${\mathbb{P}}_x$-a.s. for q.e. $x\in E$, where
$M^u\in\accentset{\circ}{\mathcal{M}}^c_{f\mbox{-}\mathrm{loc}}$
and $N^u\in{\mathcal{N}}_{c,f\mbox{-}\mathrm{loc}}$.
If, further, $u$ can be extended to be a real-valued function on
$E_{\partial}$ [without assuming $u(\partial)=0$], then the
decomposition (\ref{eqn:fuku}) holds for all
$t\in[0,\infty[$ ${\mathbb{P}}_x$-a.s. for q.e. $x\in E$.

In order to define the stochastic integrals by Dirichlet processes, we
have to
investigate the structure of $\accentset{\circ}{\mathcal{M}}^d$. For
this we
introduce the spaces
$\accentset{\circ}{\mathcal{J}}$, $\mathcal{J}$, $\widehat
{\mathcal{J}}$ of jump functions
\begin{eqnarray*}
\widehat{\mathcal{J}}&:=& \{\phi\dvtx E_{\partial}\times E_{\partial
}\to\mathbb{R}  \mid
\phi\mbox{ is a Borel measurable function such that}\\
&&\hspace*{64.8pt}\phi(x,x)=0\mbox{ for }x\in E_{\partial}\mbox{ and }
N(\mathbf{1}_{E\times E}\phi^2)\mu_H\in S \},\\
\mathcal{J}&:=& \{\phi\dvtx E_{\partial}\times E_{\partial}\to\mathbb
{R}  \mid
\phi\mbox{ is a Borel measurable function such that}\\
&&\hspace*{90.6pt}\phi(x,x)=0\mbox{ for }x\in E_{\partial}\mbox{ and }
N(\phi^2)\mu_H\in S \}
\end{eqnarray*}
and $\accentset{\circ}{\mathcal{J}}:=
\{\phi\in\mathcal{J}\mid
\int_EN(\phi^2)\,d\mu_H<\infty\}$.
Clearly $\accentset{\circ}{\mathcal{J}} \subset\mathcal{J}
\subset\widehat{\mathcal{J}}$, and
for $\phi\in\widehat{\mathcal{J}}$, we see
$\mathbf{1}_{E\times E}\phi\in\mathcal{J}$.
Further we set $\mathcal{J}_{\mathrm{as}}:=\{\phi\in\mathcal
{J}\mid\widetilde{\phi}=0$
$J\mbox{-a.e. on }E\times E\}$ and
$\mathcal{J}_*:=\{ \phi\in\mathcal{J}\mid N(\mathbf{1}_{E\times
E}|\overline{\phi
}|^2)\mu_H\in S\}$,
$\accentset{\circ}{\mathcal{J}}_{ \mathrm{as}}=\accentset{\circ}{\mathcal{J}}\cap
  \mathcal{J}_{\mathrm{as}}$ and
$\accentset{\circ}{\mathcal{J}}_*=\accentset{\circ}{\mathcal{J}}\cap  \mathcal{J}_*$.
Here $\overline\phi(x,y):=\phi(y,x)$
for $x,y\in E_{\partial}$, $\widetilde\phi:=(\phi+\overline{\phi
})/2$ on
$E_{\partial}\times E_{\partial}$. Clearly, ${\mathcal{J}}_{\mathrm
{as}}\subset\mathcal{J}_*$ and
$\accentset{\circ}{\mathcal{J}}_{ \mathrm{as}}\subset\accentset{\circ}{\mathcal{J}}_*$.
Similarly, we can define $\widehat{\mathcal{J}}_{\mathrm{as}}$ and
$\widehat{\mathcal{J}}_*$
by replacing $\mathcal{J}$ with
$\widehat{\mathcal{J}}$ in its definitions.
Moreover, for $\phi\in{\mathcal{J}}_*$ (resp., $\phi\in\widehat
{\mathcal{J}}_*$), we see
$\mathbf{1}_{E\times E}\overline{\phi}\in\mathcal{J}_*$
(resp., $\overline{\phi
}\in\widehat{\mathcal{J}}_*$).
For\vspace*{1pt} $\phi\in\widehat{\mathcal{J}}$ and $\ell\in\mathbb{N}$, we
write $\phi
_{\ell}:=\phi\mathbf{1}_{\{|\phi|>1/\ell\}}$.
For $\phi,\psi\in\mathcal{J}$ (resp., $\phi,\psi\in\widehat
{\mathcal{J}}$),
we write $\phi\sim\psi$ if $\phi=\psi$ $J^*$-a.e. on
$E\times E_{\partial}$ (resp., $\phi=\psi$ $J$-a.e. on $E\times E$),
where $J^*$ is the measure on $E\times E_{\partial}$ defined by
$J^*(dx\,dy):=\frac12 N(x,dy)\mu_H(dx)$. Then $\sim$ is an equivalence
relation and denote by
$\accentset{\circ}{\mathcal{J}}/\sim$, $\mathcal{J}/\sim$,
$\widehat{\mathcal{J}}/\sim$ the
families of equivalence classes.
\begin{lem}\label{lem:purediscontimartincorrespondince}
There exists a one-to-one correspondence between $\accentset{\circ}{\mathcal{J}}
/\sim$ and
$\accentset{\circ}{\mathcal{M}}^d$ which is characterized by the
relation that
for $\phi\in\accentset{\circ}{\mathcal{J}}$ (resp., $M\in\accentset{\circ}{\mathcal{M}}^d$),
there exists $M\in\accentset{\circ}{\mathcal{M}}^d$ (resp., $\phi\in
\accentset{\circ}{\mathcal{J}}$) such that
$\mathbf{e}(M)=\frac12 \int_EN(\phi^2)(x)\mu_H(dx)$
and
$M_t-M_{t-}=\phi(X_{t-},X_t)$ for all $t\in[0,\infty[$ ${\mathbb{P}}
_x$-a.s. for q.e. $x\in E$.
Moreover, $\langle M\rangle_t=\int_0^tN(\phi^2)(X_s)\,dH_s$ for all
$t\in[0,\infty[$ ${\mathbb{P}}_x$-a.s. for q.e. $x\in E$.
\end{lem}
\begin{pf}
Take $\phi\in\accentset{\circ}{\mathcal{J}}$ and set
\[
M_t^{\ell}:=\sum_{0<s\leq t}\phi_{\ell}(X_{s-},X_s)-\int_0^t
\int_{E_{\partial}}\phi_{\ell}(X_s,y)N(X_s,dy)\,dH_s.
\]
Then we can obtain $M^{\ell}\in\accentset{\circ}{\mathcal{M}}^d$
and
%
\begin{eqnarray}
\label{eq:qudraesti}
\langle M^{\ell}\rangle_t&=&\int_0^t\int_{E_{\partial}}\phi_{\ell}^2(X_s,y)
N(X_s,dy)\,dH_s,\\
\label{eq:energyesti}
\mathbf{e}(M^{\ell})&=&\frac12 \int_E\int_{E_{\partial}}\phi_{\ell
}^2(x,y)N(x,dy)\mu_H(dx).
\end{eqnarray}
Indeed, we easily see that $M^{\ell}$ is an MAF.
If we let $T_n^{\ell}:=\inf\{t>0\mid|M_t^{\ell}|\geq n\}$, then $\{
T_n^{\ell}\}$
is an increasing sequence of $(\mathscr{F}_t)$-stopping times such
that $\lim_{n\to\infty}T_n^{\ell}=\infty$.
Then we have $|M_{t\wedge T_n^{\ell}}^{\ell}|\leq|M_{t\wedge
T_n^{\ell}-}^{\ell}|+
|\phi(X_{t\wedge T_n^{\ell}-}, X_{t\wedge T_n^{\ell}})|\leq n+
|\phi(X_{t\wedge T_n^{\ell}-}$, $ X_{t\wedge T_n^{\ell}})|$, which
implies that for $\nu\in S_{00}$
\begin{eqnarray*}
{\mathbb{E}}_{\nu}[(M_{t\wedge T_n^{\ell}}^{\ell})^2]&\leq&2n^2\nu(E)+
2{\mathbb{E}}_{\nu} \biggl[\sum_{s\leq t}\phi^2(X_{s-},X_s) \biggr] \\
&=& 2n^2\nu(E)+ 2{\mathbb{E}}_{\nu} \biggl[\int_0^t\int_{E_{\partial
}}\phi^2(X_s,y)
N(X_s,dy)\,dH_s \biggr] \\
&\leq&2n^2\nu(E)+2(1+t)\|U_1\nu\|_{\infty}\int_EN(\phi^2)(x)\mu_H(dx)
<\infty.
\end{eqnarray*}
That is, $t\mapsto M_{t\wedge T_n^{\ell}}^{\ell}$ is a square
integrable purely discontinuous ${\mathbb{P}}_{\nu}$-martingale for
each $n$.
By Corollary A.3.1 in \cite{FOT},
\[
(M_{t\wedge T_n^{\ell}}^{\ell})^2-\sum_{s\leq t}(\Delta M_{s\wedge
T_n^{\ell}}^{\ell})^2
= (M_{t\wedge T_n^{\ell}}^{\ell})^2-\sum_{s\leq t\wedge T_n^{\ell
}}\phi_{\ell}^2(X_{s-},X_s)
\]
is a ${\mathbb{P}}_{\nu}$-martingale (also a ${\mathbb
{P}}_x$-martingale for
q.e. $x\in E$), which yields that
\begin{eqnarray*}
{\mathbb{E}}_{\nu}[(M_t^{\ell})^2]&\leq&\mathop{\underline{\lim
}}_{n\to\infty} {\mathbb{E}}_{\nu
}[(M_{t\wedge T_n^{\ell}}^{\ell})^2]
=\mathop{\underline{\lim}}_{n\to\infty}{\mathbb{E}}_{\nu} \biggl[\sum
_{s\leq t\wedge T_n^{\ell}}
\phi_{\ell}^2(X_{s-},X_s) \biggr] \\
&=&{\mathbb{E}}_{\nu} \biggl[\sum_{s\leq t}\phi_{\ell}^2(X_{s-},X_s) \biggr]
\leq{\mathbb{E}}
_{\nu} \biggl[\int_0^t\int_{E_{\partial}}\phi^2(X_s,y)
N(X_s,dy)\,dH_s \biggr] \\
&\leq& (1+t)\|U_1\nu\|_{\infty}\int_E\int_{E_{\partial}} \phi
^2(x,y)N(x,dy)\mu_H(dx)<\infty.
\end{eqnarray*}
Thus, $M_t^{\ell}$ is a square integrable MAF. Since
$\{M_{t\wedge T_n^{\ell}}^{\ell}\}_{n=1}^{\infty}$ is $L^2({\mathbb
{P}}_{\nu
})$-bounded,
by use of the Banach--Saks theorem, we have the equality
\[
{\mathbb{E}}_{\nu}[(M_t^{\ell})^2]=
{\mathbb{E}}_{\nu} \biggl[\int_0^t\int_{E_{\partial}}\phi_{\ell}^2(X_s,y)
N(X_s,dy)\,dH_s \biggr]
\]
for all $\nu\in S_{00}$.
We then have the same equation for q.e. $x\in E$ by replacing $\nu$
with $x$.
Hence $M^{\ell}\in\accentset{\circ}{\mathcal{M}}$,
(\ref{eq:qudraesti}) and (\ref{eq:energyesti}). Note that
there exists a sequence $\{T_n\}$ of totally inaccessible times such that
$\{(t,\omega)\mid M_t^{\ell}-M_{t-}^{\ell}\ne0\}=\bigcup
_{n=1}^{\infty}\auf  T_n\zu $. This yields
that $M^{\ell}=M^{\ell,d}\in\accentset{\circ}{\mathcal{M}}^d$ in
view of
Theorem A.3.9 in \cite{FOT}.
Moreover, we see that $\{M^{\ell}\}_{\ell=1}^{\infty}$ is an
\textbf{e}-Cauchy sequence in
$\accentset{\circ}{\mathcal{M}}^d$. Denote by $M\in\accentset{\circ}{\mathcal{M}}^d$
its limit. Then there exists a subsequence
$\{\ell_k\}$ such that $M^{\ell_k}$ converges to $M$ uniformly on
each compact subinterval of
$[0,\infty[$ ${\mathbb{P}}_x$-a.s. for q.e. $x\in E$.
We see for each~$\ell$, $M_t^{\ell}-M_{t-}^{\ell}=
\phi_{\ell}(X_{t-},X_t)$ for all $t\in[0,\infty[$ ${\mathbb{P}}_x$-a.s. for
q.e. $x\in E$.
Therefore we have the desired result.
Conversely take an $M\in\accentset{\circ}{\mathcal{M}}^d$. Then, by
Lemma 3.2 in
\cite{cftyz},
there exists a Borel function
$\phi$ defined on $E_{\partial}\times E_{\partial}$
with $\phi(x,x)=0$ for $x\in E_{\partial}$
such that $M_t-M_{t-}=\phi(X_{t-},X_t)$, $t\in[0,\zeta_p[$ ${\mathbb
{P}}_m$-a.e.,
where $\zeta_p$ is the predictable part of $\zeta$ defined by
$\zeta_p:=\zeta$ if $\zeta<\infty$ and $X_{\zeta-}=\partial$, and
$\zeta_p:=\infty$ otherwise.
For $L=f*M^u$ with $f\in C_0(E)$, $u\in\mathcal{F}\cap C_0(E)$,
we have $\Delta L_{\zeta_p}=f(X_{\zeta_p-})\Delta M_{\zeta_p}^u=0$
for $\zeta_p<\infty$.
In view of Lemma 5.6.3 and Theorem 5.2.1 in \cite{FOT},
we see $\Delta M_{\zeta_p}=0$ for $\zeta_p<\infty$, which implies
$M_t-M_{t-}=\phi(X_{t-},X_t)$, $t\in[0,\infty[$ ${\mathbb{P}}_m$-a.e.
From this, we have
\begin{eqnarray*}
\int_E\int_{E_{\partial}}\phi^2(x,y)N(x,dy)\mu_H(dx)&=&
\lim_{t\to0}\frac1{t}{\mathbb{E}}_m \biggl[\int_0^t\int_{E_{\partial
}}\phi
^2(X_s,y)N(X_s,dy)\,dH_s \biggr]\\
&=&\lim_{t\to0}\frac1{t}{\mathbb{E}}_m \biggl[\sum_{s\leq t}\phi
^2(X_{s-},X_s) \biggr]\\
&=&\lim_{t\to0}\frac1{t}{\mathbb{E}}_m \biggl[\sum_{s\leq t}(M_s-M_{s-})^2
\biggr]\\
&=& \lim_{t\to0}\frac1{t}{\mathbb{E}}_m[M_t^2]=
2\mathbf{e}(M)<\infty,
\end{eqnarray*}
where we use Corollary A.3.1 in \cite{FOT}.
Going back to the first argument, we can construct
$\widetilde{M}\in\accentset{\circ}{\mathcal{M}}^d$ such that
$\widetilde{M}_t-\widetilde{M}_{t-}=\phi(X_{t-},X_t)$, $t\in
[0,\infty[$ ${\mathbb{P}}
_x$-a.s. for q.e. $x\in E$.
Applying Corollary A.3.1 in \cite{FOT} to $M-\widetilde{M}\in
\accentset{\circ}{\mathcal{M}}^d$,
we obtain
\[
\mathbf{e}(M-\widetilde{M})=\lim_{t\to0}\frac1{t}{\mathbb{E}}_m
[(M_t-\widetilde
{M}_t)^2 ]=
\lim_{t\to0}\frac1{t}{\mathbb{E}}_m \biggl[\sum_{s\leq t} \bigl(\Delta
(M-\widetilde{M})_s
\bigr)^2 \biggr]=0,
\]
which implies the converse assertion.
\end{pf}
\begin{cor}\label{cor:MAFreversion}
Take $\phi\in\accentset{\circ}{\mathcal{J}}$ and set $\overline
{\phi
}(x,y):=\phi(y,x)$ for
$x,y\in E_{\partial}$.
Then $\mathbf{1}_{E\times E}\overline{\phi}\in\accentset{\circ}{\mathcal{J}}$,
in particular, there exists $K\in\accentset{\circ}{\mathcal{M}}^j$
such that
$K_t-K_{t-}=-\mathbf{1}_{E\times E}(\phi+\overline{\phi})(X_{t-},X_t)$
$t\in\ ]0,\infty[$ ${\mathbb{P}}_x$-a.s. for q.e. $x\in E$.
\end{cor}
\begin{pf}
The assertion is clear from
\begin{eqnarray*}
\int_E\int_E\overline{\phi}{}^2(x,y)N(x,dy)\mu_H(dx)&=&\int_E\int
_E{\phi}^2(x,y)N(x,dy)\mu_H(dx)\\
&\leq&\int_{E}N({\phi}^2)(x)\mu_H(dx)<\infty.
\end{eqnarray*}
\upqed\end{pf}

From this corollary, we have $\accentset{\circ}{\mathcal{J}}\subset
\mathcal{J}_*$.
\begin{lem}\label{lem:jumplocalonG}
Take a Borel function $\phi\dvtx E_{\partial}\times E_{\partial}\to
\mathbb{R}$
with $\phi(x,x)=0$ for $x\in E_{\partial}$.
The following are equivalent under $\phi\in\mathcal{J}$:
\begin{enumerate}[(5)]
\item[(1)]
$\phi(X_{t-},X_t)=0$ for all $t\leq
\sigma_{E\setminus G}$
${\mathbb{P}}_x$-a.s. for q.e. $x\in E$.
\item[(2)]
$\phi(X_{t-},X_t)=0$ for all $t\leq
\sigma_{E\setminus G}$
${\mathbb{P}}_m$-a.e.
\item[(3)]
$\int_0^tN(\phi^2)(X_s)\,dH_s=0$ for all
$t<\sigma_{E\setminus G}$
${\mathbb{P}}_x$-a.s. for q.e. $x\in E$.
\item[(4)]
$\int_0^tN(\phi^2)(X_s)\,dH_s=0$ for all
$t<\sigma_{E\setminus G}$
${\mathbb{P}}_m$-a.e.
\item[(5)]
$\phi=0$ $J^*$-a.e. on $G\times
E_{\partial}$.
\end{enumerate}
Replacing $\sigma_{E\setminus G}$ [resp., $\phi$ in
(3), (4)]
with $\tau_G$ (resp., $\mathbf{1}_{E\times E}\phi$), we have a
similar equivalence under
$\mathbf{1}_{E\times E}\phi\in{\mathcal{J}}$, where
the last condition is that $\phi=0$ $J^*$-a.e. on $G\times E$.
\end{lem}
\begin{pf}
The implication
(1) $\Longrightarrow$ (2) is
trivial and
(3) $\Longleftrightarrow$ (4)
follows from
Lemma 5.1.10(iii) in \cite{FOT}.
We first show
(2) $\Longrightarrow$ (3).
Suppose
$\phi(X_{t-},X_t)=0$ for all $t\in\ ]0,\sigma_{E\setminus G}]$
${\mathbb{P}}_m$-a.e.
Then we see
$\phi(X_{\sigma_{E\setminus G}-},X_{\sigma_{E\setminus G}})=0$
${\mathbb{P}}
_m$-a.e. on
$\{\sigma_{X\setminus G}<\infty\}$.
So
$\phi(X_{t\wedge\sigma_{E\setminus G}-},X_{t\wedge\sigma
_{E\setminus G}})=0$ for all
$t\in\ ]0,\infty[$
${\mathbb{P}}_m$-a.e. From the property of L\'evy system (see Appendix
(A) in
\cite{CKmixed} or
the formula with $Y_t=1_{]0,T]}(t)$ at line -9 on page 346 in
\cite{Shargene}),
we have for each $t>0$
\[
{\mathbb{E}}_m \biggl[\int_0^{t\wedge\sigma_{E\setminus G}}
N(\phi^2)(X_s)\,dH_s \biggr]={\mathbb{E}}_m \biggl[\sum_{s\leq t\wedge\sigma
_{E\setminus
G}}\phi^2(X_{s-},X_s) \biggr]=0,
\]
which implies (4), hence (3).
(3) also yields
$\int_GN(\phi^2)\,d\mu_H=0$ by Lem\-ma~5.1.10(iii) in \cite{FOT}, and
in particular, we obtain (5).
Conversely suppose (5), that is, $\int_GN(\phi
^2)\,d\mu_H=0$.
Then, we can obtain (1)
by way of the inequality (\ref{eq:fineversionoflem5.1.9}) and
the property of L\'evy system used above.
\end{pf}
\begin{cor}\label{cor:MAFjkchara}
Take an MAF $M\in\accentset{\circ}{\mathcal{M}}^d$ and the
associated $\phi
\in\accentset{\circ}{\mathcal{J}}$.
Set $\phi_{\partial}(x,y):=\phi(x,y)\mathbf{1}_{\{\partial\}}(y)$.
Then the following are equivalent:
\begin{enumerate}[(4)]
\item[(1)] 
$M\in\accentset{\circ}{\mathcal{M}}^j$.
\item[(2)] 
$\phi(x,\partial)=0$
$\kappa$-a.e. $x\in E$.
\item[(3)] 
$\int_0^{\cdot}N(\phi
_{\partial}^2)(X_s)\,dH_s\equiv0$
${\mathbb{P}}_x$-a.s. for q.e. $x\in E$.
\item[(4)] 
$\int_0^{\cdot}N(\phi
_{\partial}^2)(X_s)\,dH_s\equiv0$
${\mathbb{P}}_m$-a.e.
\end{enumerate}
Set $\phi_{E}(x,y):=\phi(x,y)\mathbf{1}_{E}(y)$.
Then the following are equivalent:
\begin{enumerate}[(4$^*$)]
\item[(1$^*$)]
$M\in \accentset{\circ}{\mathcal{M}}^{\kappa}$.
\item[(2$^*$)]
$\phi (x,y)=0$ $J$-a.e. $(x,y)\in E\times E$.
\item[(3$^*$)]
$\int_0^{\cdot}
N(\phi_{E}^2)(X_s)\,dH_s\equiv0$
${\mathbb{P}}_x$-a.s. for q.e. $x\in E$.
\item[(4$^*$)]
$\int_0^{\cdot}
N(\phi_{E}^2)(X_s)\,dH_s\equiv0$
${\mathbb{P}}_m$-a.e.
\end{enumerate}
\end{cor}
\begin{pf}
(1) $\Longleftrightarrow$ (2)
is clear from
$\mathbf{e}(M,M^{u,\kappa})=-\frac{1}{2}
\int_E\phi(x,\partial)u(x)\kappa(dx)$ for $u\in\mathcal{F}_e$.
Here we use the fact that $\mathcal{F}$ is dense in $L^2(E;\kappa)$.
(1$^*$) $\Longleftrightarrow$ (2$^*$)
is clear from
(1) $\Longleftrightarrow$ (2) and
$\mathbf{e}(M,N)=\frac12
\int_E\int_{E_{\partial}}\phi(x,y)\psi(x,y)N(x,dy)\mu_H(dx)$.
The rest implications hold true for general $\phi\in\mathcal{J}$ and are
clear in view of the uniqueness of the Revuz correspondence and
Lemma 5.1.10(iii) in \cite{FOT}.
\end{pf}

Let $\mathcal{M}_{\mathrm{loc}}$ be the space of locally square integrable MAFs and
$\mathcal{M}_{\mathrm{loc}}^{\auf 0,\zeta\auf }$ the space of locally square integrable
MAFs on $\auf 0,\zeta\auf $.
That is, $M\in\mathcal{M}_{\mathrm{loc}}$ (resp.,  $M\in\mathcal{M}_{\mathrm{loc}}^{\auf 0,\zeta\auf }$)
if and only if there exists an increasing
sequence $\{T_n\}$ (resp., $\{S_n\}$) of $(\mathscr{F}_t)$-stopping
times and $\{M^{(n)}\}\subset\mathcal{M}$
such that $\lim_{n\to\infty}T_n=\infty$ (resp., $\lim_{n\to\infty
}S_n=\zeta$) and for each $n\in\mathbb{N}$,
$M_{t\wedge{T_n}}=M^{(n)}_{t\wedge T_n}$
(resp.,\break $M_{t\wedge{S_n}}\mathbf{1}_{\{
t\wedge S_n<\zeta\}}=M^{(n)}_{t\wedge S_n}\mathbf{1}_{\{t\wedge
S_n<\zeta\}}$) for all
$t\in[0,\infty[$ ${\mathbb{P}}_x$-a.s. for q.e. $x\in E$.
Let $\mathcal{M}^{c}_{\mathrm{loc}}$ (resp., $\mathcal{M}^{d}_{\mathrm{loc}}$) be the
space of locally
square integrable continuous
(resp., purely discontinuous) MAFs. That is, for $M\in\mathcal{M}^{c}_{\mathrm{loc}}$
(resp., $M\in\mathcal{M}^{d}_{\mathrm{loc}}$),
we can take $\{M^{(n)}\}$ from $\mathcal{M}^c$ (resp., $\mathcal
{M}^d$) in the above definition.
Similarly, we can define the space
$\mathcal{M}^{c,\auf 0,\zeta\auf }_{\mathrm{loc}}$ (resp., $\mathcal{M}^{d,\auf 0,\zeta\auf }_{\mathrm{loc}}$) of
locally square integrable continuous
(resp., purely discontinuous) MAFs on $\auf 0,\zeta\auf $.
For every $M\in\mathcal{M}_{\mathrm{loc}}^{\auf 0,\zeta\auf }$, its quadratic
variational process $\langle M\rangle$ can be defined to be a PCAF
(Proposition 2.8 in \cite{CFKZStoch}), and $M$
is decomposed to $M=M^c+M^d$ (Theorem 8.23 in \cite{HWY}),
where $M^c\in\mathcal{M}^{c,\auf 0,\zeta\auf }_{\mathrm{loc}}$,
$M^d\in\mathcal{M}^{d,\auf 0,\zeta\auf }_{\mathrm{loc}}$ have the property
$\langle
M^c,M^d\rangle\equiv0$.
The next theorem is a natural extension of Lemma \ref
{lem:purediscontimartincorrespondince}.
\begin{theorem}\label{thm:purediscontimartincorrespondince}
There exists a one-to-one correspondence between $\mathcal{J}/\sim$
(resp., $\widehat{\mathcal{J}}/\sim$) and
$\mathcal{M}^d_{\mathrm{loc}}$ (resp., $\mathcal{M}^{d,\auf 0,\zeta\auf }_{\mathrm{loc}}$)
which is characterized by the relation that
for $\phi\in{\mathcal{J}}$ (resp., $\widehat{\mathcal{J}}$)
there exists $M\in\mathcal{M}^d_{\mathrm{loc}}$ (resp., $\mathcal{M}^{d,\auf 0,\zeta\auf }_{\mathrm{loc}}$)
such that
$M_t-M_{t-}=\phi(X_{t-},X_t)$ for all $t\in[0,\infty[$ (resp., $t\in
[0,\zeta[$)
${\mathbb{P}}_x$-a.s. for q.e. $x\in E$. Conversely
for each $M\in\mathcal{M}^d_{\mathrm{loc}}$ (resp., $\mathcal{M}^{d,\auf 0,\zeta\auf }_{\mathrm{loc}}$),
there exists a
$\phi\in{\mathcal{J}}$ (resp., $\widehat{\mathcal{J}}$) such that
the same equation holds.
Moreover, we have $\langle M\rangle_t=\int_0^t\int_{E_{\partial
}}\phi^2(X_s,y)N(X_s,dy)\,dH_s$ for all
$t\in[0,\infty[$ ${\mathbb{P}}_x$-a.s. for q.e. $x\in E$.
\end{theorem}
\begin{pf}
We only prove the correspondence between $\mathcal{J}/\sim$ and
$\mathcal{M}^d_{\mathrm{loc}}$.
The proof of the correspondence between $\widehat{\mathcal{J}}/\sim
$ and $\mathcal{M}
^{d,\auf 0,\zeta\auf }_{\mathrm{loc}}$
is similar
by replacing $\sigma_{E\setminus F_k}$ with $\tau_{F_k}$.
Suppose $\phi\in\mathcal{J}$.
Take an $\mathcal{E}$-nest $\{F_k\}$ of closed sets
such that
$\mathbf{1}_{F_k}N(|\phi|^2
)\mu_H\in S_{00}$. Then $\mathbf{1}_{F_k}\phi\in\accentset{\circ}{\mathcal{J}}$ for each $k\in\mathbb{N}$,
where $(\mathbf{1}_{F_k}\phi)(x,y):=\mathbf{1}_{F_k}(x)\phi(x,y)$,
$x,y\in E_{\partial}$, and
there exists an $M^{(k)}\in\accentset{\circ}{\mathcal{M}}^d$ such that
$M_t^{(k)}-M_{t-}^{(k)}=\mathbf{1}_{F_k}(X_{t-})\phi(X_{t-},X_t)$ for
all $t\in[0,\infty[$
${\mathbb{P}}_x$-a.s. for q.e. $x\in E$.
Such $M^{(k)}$ is an $\mathbf{e}$-convergent limit of $\{M^{(k),\ell
}\}_{\ell=1}^{\infty}$, where
\[
M_t^{(k),\ell}:=\sum_{0<s\leq t}\mathbf{1}_{F_k}(X_{s-})\phi_{\ell
}(X_{s-},X_s)
-\int_0^t\mathbf{1}_{F_k}(X_s)N(\phi_{\ell})(X_s)
\,dH_s.
\]
This yields that for $j>k$, $M_t^{(j)}=M_t^{(k)}$ for $t<\sigma
_{E\setminus F_k}$, more strongly\break
$M^{(j)}_{t\wedge\sigma_{E\setminus F_k}}=M^{(k)}_{t\wedge\sigma
_{E\setminus F_k}}$ because of
$X_{\sigma_{E\setminus F_k}-}\in F_k\cup\{\partial\}$.
Hence $M$ defined by $M_t:=M_t^{(k)}$ for $t<\sigma_{E\setminus F_k}$ satisfies
$M_{t\wedge\sigma_{E\setminus F_k}}=M^{(k)}_{t\wedge\sigma
_{E\setminus F_k}}$, which
implies $M\in\mathcal{M}^d_{\mathrm{loc}}$, because $\{F_k\}$ is also a strict
$\mathcal{E}
$-nest.

Conversely suppose $M\in\mathcal{M}^d_{\mathrm{loc}}$.
Then there exists a sequence $\{M^{(n)}\}$
of square integrable purely discontinuous MAFs
and an increasing sequence $\{T_n\}$ of stopping
times such that $\lim_{n\to\infty}T_n=\infty$ and
$M_{t\wedge T_n}=M^{(n)}_{t\wedge T_n}$ for all $t\in[0,\infty[$
${\mathbb{P}}_x$-a.s. for q.e. $x\in E$.
By an argument in the proof of Proposition 2.8 in \cite{CFKZStoch}, we
can construct a quadratic
variational process $\langle M\rangle$, which is a PCAF,
and a nest
$\{F_k\}$ of closed sets such that $\mathbf{1}_{F_k\cup\{\partial\}
}*M\in\accentset{\circ}{\mathcal{M}}$
and $\mathbf{e}(\mathbf{1}_{F_k\cup\{\partial\}}*M)={\frac
{1}{2}}\mu_{\langle M\rangle}(F_k)$.
Note that $\mathbf{1}_{\{\partial\}}*M=0$ because\break $\int_0^t\mathbf
{1}_{\{\partial\}}(X_s)\,d\langle M\rangle_s=0$.
We remark that $\langle M,N\rangle\equiv0$ for all $N\in\mathcal
{M}^{c}_{\mathrm{loc}}$,
which implies
$\mathbf{1}_{F_k\cup\{\partial\}}*M\in\accentset{\circ}{\mathcal{M}}^d$, hence
$M\in\accentset{\circ}{\mathcal{M}}^d_{f\mbox{-}\mathrm{loc}}$.
As in the proof of Proposition 2.8 in \cite{CFKZStoch}, we see
$M_t=(\mathbf{1}_{F_k\cup\{\partial\}}*M)_t$ for $t\leq\sigma
_{E\setminus F_k}$.
Indeed, we have this from the assertion for $t< \sigma_{E\setminus
F_k}$ and
$\Delta M_{\sigma_{E\setminus F_k}}=\mathbf{1}_{F_k\cup\{\partial\}
}(X_{\sigma_{E\setminus F_n}-})
\Delta M_{\sigma_{E\setminus F_k}}=\Delta(\mathbf{1}_{F_k\cup\{
\partial\}}*M)_{\sigma_{E\setminus F_k}}$.
By Lemma \ref{lem:purediscontimartincorrespondince},
there exists a Borel function $\phi_k\in\accentset{\circ}{\mathcal{J}}$ such that
$\mathbf{1}_{F_k\cup\{\partial\}}(X_{t-})(M_t-M_{t-})=\phi
_k(X_{t-},X_t)$ for $t\in[0,\infty[$
${\mathbb{P}}_x$-a.s. for q.e. $x\in E$. From this, for $j>k$, we see that
$\phi_k(X_{t-},X_t)=\phi_j(X_{t-},X_t)$ for
$t\leq\sigma_{E\setminus F_k}$ ${\mathbb{P}}_x$-a.s. for q.e. $x\in E$.
Let $G_k$ be the fine interior of $F_k$.
By Lemma \ref{lem:jumplocalonG}, for $j>k$ we have
$\phi_k=\phi_j$ $J^*$-a.e. on $G_k\times E_{\partial}$.
So we can define $\phi$ on $E_{\partial}\times E_{\partial}$ such that
$\phi=\phi_k$ $J^*$-a.e. on $G_k\times E_{\partial}$.
From Lemma \ref{lem:smoothfinelyopen}, we see $N(\phi^2)\mu_H\in S$.
Applying Lemma \ref{lem:jumplocalonG} again,
$\phi(X_{t-},X_t)=\phi_k(X_{t-},X_t)=M_t-M_{t-}$
for all $t\leq\sigma_{E\setminus G_k}$.
Moreover, we see $\langle\mathbf{1}_{F_k\cup\{\partial\}}*M\rangle_t=
\int_0^t\int_{E_{\partial}}
\phi_k^2(X_s,y)N(X_s,dy)\,dH_s$.
Therefore we obtain the desired assertion.
\end{pf}
\begin{cor}\label{cor:MAFreversion*}
For $\phi\in\widehat{\mathcal{J}}_*$, there exists a
$K\in\mathcal{M}^{d}_{\mathrm{loc}}$ such that $K_t-K_{t-}=-\mathbf
{1}_{E\times
E}(\phi+\overline{\phi})(X_{t-},X_t)$ for all $t\in\ ]0,\infty[$
${\mathbb{P}}
_x$-a.s. for q.e. $x\in E$.
\end{cor}
\begin{pf}
This is clear from that $\phi\in\widehat{\mathcal{J}}_*$
implies $\mathbf{1}_{E\times E}\phi, \mathbf{1}_{E\times E}\overline
{\phi}\in\mathcal{J}_*$.
\end{pf}
\begin{rem}\label{rem:coincidenceMAFs}
A similar argument of the proof of Theorem \ref
{thm:purediscontimartincorrespondince} yields
\begin{eqnarray*}
\mathcal{M}^{d}_{\mathrm{loc}}&=&\bigl\{M\mid\mbox{ there exists }\{G_n\}\in\Theta
\mbox{
and }M^{(n)}\in\accentset{\circ}{\mathcal{M}}^d\mbox{ such that }
\\
&&\hspace*{5.1pt} M_t=M_t^{(n)}\mbox{ for all }t\leq\sigma
_{E\setminus G_n}
\mbox{ and }n\in\mathbb{N}, {\mathbb{P}}_x\mbox{-a.s. for
q.e. }x\in E\bigr\}, \\
\mathcal{M}^{d,\auf 0,\zeta\auf }_{\mathrm{loc}}&=&\bigl\{M\mid\mbox{ there exists
}\{G_n\}
\in\Theta\mbox{ and }M^{(n)}\in\accentset{\circ}{\mathcal{M}}^d\mbox{
such that } \\
&&\hspace*{5.1pt} M_t=M_t^{(n)}\mbox{ for all }t\leq\tau_{G_n}
\mbox{ and }n\in\mathbb{N}, {\mathbb{P}}_x\mbox{-a.s. for
q.e. }x\in E\bigr\}.
\end{eqnarray*}
These show $\mathcal{M}^{d}_{\mathrm{loc}}\subset\accentset{\circ}{\mathcal{M}}
^d_{f\mbox{-}\mathrm{loc}}$ and
$\mathcal{M}^{d,\auf 0,\zeta\auf }_{\mathrm{loc}}\subset\accentset{\circ}{\mathcal{M}}^d_{\mathrm{loc}}$.
We also have the coincidences $\mathcal{M}^{c}_{\mathrm{loc}}=\accentset{\circ}{\mathcal{M}}^c_{f\mbox{-}\mathrm{loc}}$ and $\mathcal
{M}^{c,\auf 0,\zeta\auf }_{\mathrm{loc}}=
\accentset{\circ}{\mathcal{M}}^c_{\mathrm{loc}}$.
Indeed, the inclusion
$\mathcal{M}^{c}_{\mathrm{loc}}\subset\accentset{\circ}{\mathcal{M}}^c_{
f\mbox{-}\mathrm{loc}}$
can be obtained in the same way of the proof of Theorem \ref
{thm:purediscontimartincorrespondince}. The converse inclusion is
easily confirmed
from the continuity of $M\in\accentset{\circ}{\mathcal{M}}
^c_{f\mbox{-}\mathrm{loc}}$
and ${\mathbb{P}}_x(\lim_{n\to\infty}\sigma_{E\setminus
G_n}=\infty)=1$
for q.e. $x\in E$.
\end{rem}

The next corollary is needed to assure the uniqueness of the
generalized Fuku\-shima decomposition later.
\begin{cor}\label{cor:uniqeness}
We have $\mathcal{M}_{\mathrm{loc}}\cap\mathcal{N}_{c,f\mbox{-}\mathrm{loc}}=\{0\}$ and
$\mathcal{M}^{\auf 0,\zeta\auf }_{\mathrm{loc}}\cap\mathcal{N}_{c,\mathrm{loc}}=\{0\}$.
\end{cor}
\begin{pf}
We only prove $\mathcal{M}_{\mathrm{loc}}\cap\mathcal{N}_{c,f\mbox{-}\mathrm{loc}}=\{
0\}$. The
proof of
$\mathcal{M}^{\auf 0,\zeta\auf }_{\mathrm{loc}}\cap\mathcal{N}_{c,\mathrm{loc}}=\{0\}$ is
similar to this by replacing
$\sigma_{E\setminus F_k}$ with $\tau_{F_k}$.
Take $M\in\mathcal{M}_{\mathrm{loc}}$, $N\in\mathcal{N}_{c,f\mbox{-}\mathrm{loc}}$ and
suppose $M_t+N_t=0$ for all
$t\in[0,\infty[$ ${\mathbb{P}}_x$-a.s. for q.e. $x\in E$.
In particular, $M$ is continuous.
Let $\phi\in\mathcal{J}$ be the jump function associated to $M^d$.
As in the
proof of
Theorem \ref{thm:purediscontimartincorrespondince},
we can construct a common
$\{G_k\}\in\Theta$, $M^{(k)}\in\accentset{\circ}{\mathcal{M}}$
with its
jump function $\phi_k\in\accentset{\circ}{\mathcal{J}}$
and $N^{(k)}\in\mathcal{N}_c$ such that
$M_t=M_t^{(k)}$, $N_t=N_t^{(k)}$ for all $t<\sigma_{E\setminus G_k}$
${\mathbb{P}}_x$-a.s. for
q.e. $x\in E$, and $\phi=\phi_k$ $J^*$-a.e. on $G_k\times E_{\partial}$.
The continuity of $M$ yields $\phi(X_{t-},X_t)=0$ for all $t\in
\ ]0,\infty[$ ${\mathbb{P}}_x$-a.s. for q.e. $x\in E$.
Then we can conclude that $\phi_k(X_{t-},X_t)=0$ for all $t\in
\ ]0,\sigma_{E\setminus G_k}]$
${\mathbb{P}}_x$-a.s. for q.e. $x\in E$ by way of Lemma
\ref{lem:jumplocalonG}.
This implies that $M_{t\wedge\sigma_{E\setminus G_k}}^{(k)}=
M_{t\wedge\sigma_{E\setminus G_k}-}^{(k)}$ ${\mathbb{P}}_x$-a.s. for
q.e. $x\in E$. On the other hand,
we see that $M_{t\wedge\sigma_{E\setminus G_k}-}^{(k)}+
N_{t\wedge\sigma_{E\setminus G_k}}^{(k)}=0$ ${\mathbb{P}}_x$-a.s. for
q.e. $x\in E$ because $M_t^{(k)}+N_t^{(k)}=0$ for all
$t\in[0,\sigma_{E\setminus G_k}[$ ${\mathbb{P}}_x$-a.s. for
q.e. $x\in E$.
Therefore we obtain
\[
M_{t\wedge\sigma_{E\setminus G_k}}^{(k)}+
N_{t\wedge\sigma_{E\setminus G_k}}^{(k)}=0  \qquad\mbox{for all }t\in
[0,\infty[ \mbox{ } {\mathbb{P}}_x\mbox{-a.s. for q.e. }
x\in E.
\]
Then we can conclude that $M_t=N_t=0$ for all $t\in[0,\infty[$
${\mathbb{P}}
_x$-a.s. for q.e. $x\in E$ in view of the
argument of the proof of the uniqueness of Fukushima decomposition as
in Theorem 5.5.1 of \cite{FOT}.
\end{pf}

We define subclasses of $\mathcal{M}^d_{\mathrm{loc}}$ as follows:
\begin{eqnarray*}
\mathcal{M}^j_{\mathrm{loc}}&:=\{M\in\mathcal{M}^d_{\mathrm{loc}} \mid\phi(\cdot
,\partial)=0\mbox{ }
\kappa\mbox{-a.e.  on }E\}, \\
\mathcal{M}^{\kappa}_{\mathrm{loc}}&:=\{M\in\mathcal{M}^d_{\mathrm{loc}} \mid\phi=0\mbox{ }
J\mbox{-a.e.
on }E\times E\}.
\end{eqnarray*}
Then we have a similar statement as in Corollary \ref{cor:MAFjkchara}.
From this, we see that
$M\in\mathcal{M}^j_{\mathrm{loc}}$, $N\in\mathcal{M}^{\kappa}_{\mathrm{loc}}$
implies $\langle M,N\rangle\equiv0$ ${\mathbb{P}}_x$-a.s. for
q.e. $x\in E$.
Every $M\in\mathcal{M}_{\mathrm{loc}}$ is decomposed to $M=M^c+M^j+M^{\kappa
}$, where
$M^c\in\mathcal{M}^c_{\mathrm{loc}}$, $M^j\in M^j_{\mathrm{loc}}$, $M^{\kappa}\in
\mathcal{M}^{\kappa
}_{\mathrm{loc}}$ have the properties
$\langle M^c,M^j\rangle\equiv\langle M^j,M^{\kappa}\rangle\equiv
\langle M^{\kappa},M^c\rangle\equiv0$.
For $M\in\mathcal{M}_{\mathrm{loc}}^{\auf 0,\zeta\auf }$ with its jump function
$\phi
\in\widehat{\mathcal{J}}$, we can consider $M^j\in\mathcal
{M}_{\mathrm{loc}}^j$ (resp., $K\in\mathcal{M}
_{\mathrm{loc}}^j$) associated to $\mathbf{1}_{E\times E}\phi\in\mathcal{J}$
[resp., $-\mathbf{1}_{E\times E}(\phi+\overline{\phi})\in\mathcal
{J}$], where
$K$ is constructed in
Corollary \ref{cor:MAFreversion*}.

We introduce the subclasses $\dot{\mathcal{F}}^{\dag}_{\mathrm{loc}}$, $\dot
{\mathcal{F}
}_{\mathrm{loc}}^{\ddag}$ of $\dot{\mathcal{F}}_{\mathrm{loc}}$ as follows:
\begin{eqnarray*}
\dot{\mathcal{F}}^{\dag}_{\mathrm{loc}}&:=&\bigl\{u\in\dot{\mathcal
{F}}_{\mathrm{loc}}\mid N\bigl(\mathbf
{1}_{E\times E}\bigl(u(\cdot)-u\bigr)^2\bigr)\mu_H\in S \bigr\}, \\
\dot{\mathcal{F}}^{\ddag}_{\mathrm{loc}}&:=&\bigl\{u\in\dot{\mathcal{F}}^{\dag
}_{\mathrm{loc}}\mid
u(\partial)\in\mathbb{R}\mbox{ and }\bigl(u(\cdot)-u(\partial
)\bigr)^2\kappa\in S
\bigr\}.
\end{eqnarray*}
Clearly, $\dot{\mathcal{F}}^{\dag}_{\mathrm{loc}}$ and
$\dot{\mathcal{F}}_{\mathrm{loc}}^{\ddag}$ are linear subspaces of $\dot
{\mathcal{F}}_{\mathrm{loc}}$,
and $\mathbf{1}_{E_{\partial}},\mathbf{1}_E\in\dot{\mathcal{F}
}_{\mathrm{loc}}^{\ddag}$. By Remark 3.9 of \cite{CFKZStoch} and $\kappa\in S$,
we see $\mathcal{F}_e\cup(\dot{\mathcal{F}}_{\mathrm{loc}})_b
\subset\dot{\mathcal{F}}^{\ddag}_{\mathrm{loc}}$ by regarding $u(\partial
)\in\mathbb{R}$
for $u\in\dot{\mathcal{F}}_{\mathrm{loc}}$.
For $u,v\in\dot{\mathcal{F}}_{\mathrm{loc}}^{\ddag}$ (resp., $u,v\in\dot
{\mathcal{F}
}_{\mathrm{loc}}^{\dag}$), we see
$uv\in\dot{\mathcal{F}}_{\mathrm{loc}}^{\ddag}$ (resp., $uv\in\dot{\mathcal
{F}}_{\mathrm{loc}}^{\dag
}$) provided $u$ or $v$ is bounded.
From Theorem \ref{thm:purediscontimartincorrespondince}, for $u\in
\dot{\mathcal{F}}^{\ddag}_{\mathrm{loc}}$
(resp., $u\in\dot{\mathcal{F}}_{\mathrm{loc}}^{\dag}$),
there exists a $M^{u,d}\in\mathcal{M}^{d}_{\mathrm{loc}}(\subset
\accentset{\circ}{\mathcal{M}}^d_{f\mbox{-}\mathrm{loc}})$
(resp., $M^{u,d}\in\mathcal{M}^{d,\auf 0,\zeta\auf }_{\mathrm{loc}}$) such that
$\Delta M^{u,d}_t=\Delta u(X_t)$ for all $t\in[0,\infty[$
(resp., $[0,\zeta[$) ${\mathbb{P}}_x$-a.s. for
q.e. $x\in E$.

Moreover, we define
\begin{eqnarray*}
{\mathcal{F}}^{\dag}_{\mathrm{loc}}&:=& \biggl\{u\in{\mathcal{F}}_{\mathrm{loc}}  \Bigm|
\forall K\in\mathcal{K}, \int
_{K\times E}\bigl(u(y)-u(x)\bigr)^2J(dx\,dy)<\infty\biggr\}, \\
{\mathcal{F}}^{\ddag}_{\mathrm{loc}}&:=& \biggl\{u\in{\mathcal{F}}_{\mathrm{loc}}^{\dag}
\Bigm|  u(\partial
)\in\mathbb{R}\mbox{ and }\forall K\in\mathcal{K}, \int
_K\bigl(u(x)-u(\partial
)\bigr)^2\kappa(dx)<\infty\biggr\}.
\end{eqnarray*}
Here $\mathcal{K}$ denotes the family of all compact sets and
$\mathcal{F}_{\mathrm{loc}}$ is the space of functions locally in $\mathcal
{F}$ in the ordinary
sense (see \cite{FOT}).
Clearly, $\mathcal{F}_{\mathrm{loc}}^{\dag}\subset\dot{\mathcal
{F}}_{\mathrm{loc}}^{\dag}$ and
$\mathcal{F}_{\mathrm{loc}}^{\ddag}\subset\dot{\mathcal{F}}_{\mathrm{loc}}^{\ddag}$.
For $u\in{\mathcal{F}}_{\mathrm{loc}}$, $u\in\mathcal{F}_{\mathrm{loc}}^{\dag}$ if
and only if that
for any compact set $K$ with its relatively
compact open neighborhood $G$
\[
\int_{K\times G^c}\bigl(u(y)-u(x)\bigr)^2J(dx\,dy)<\infty.
\]
We see $\mathcal{F}_e\cup(\mathcal{F}_{\mathrm{loc}})_b\subset\mathcal
{F}_{\mathrm{loc}}^{\ddag}$, because of
$J(K\times G^c)<\infty$ and $\kappa(K)<\infty$ (see Corollary 5.1 in
\cite{Kwfunc}), where $K$ and $G$ are noted as above.
\begin{ex}
Let $(\mathcal{E},\mathcal{F})$ be a regular Dirichlet form on
$L^2(\mathbb{R}^N)$ whose
jumping measure $J$ has an expression
$J(dx\,dy)=f(|x-y|)\,dx\,dy$ such that $f$ is a locally bounded Borel
function on $[0,\infty[$ satisfying
%
\begin{equation}\label{eq:conditionjump}
\int_c^{\infty}f(r)r^{N+1}\,dr<\infty\qquad\mbox{for some } c>0.
\end{equation}
For instance, relativistic symmetric $\alpha$-stable processes
satisfy (\ref{eq:conditionjump}) (see \cite{CKmixed}).
Then each coordinate function $\ell_k(x):=x_k$ satisfies
$\ell_k\in{\mathcal{F}}^{\dag}_{\mathrm{loc}}$ ($k=1,2,\ldots,N$)
under~(\ref{eq:conditionjump}).
Indeed, for any compact set $K$ and its relatively compact open
neighborhood $G$ with $d(K,G^c):=
{\inf_{x\in K,y\in G^c}}|x-y|>0$,
\begin{eqnarray*}
\int_{K\times G^c}|\ell_k(x)-\ell_k(y)|^2J(dx\,dy)&\leq&
\int_K\int_{G^c}|x-y|^2f(|x-y|)\,dy\,dx \\
&\leq&|K|\sigma(\mathbb{S}^{N-1})\int_{d(K,G^c)}^{\infty
}r^{N+1}f(r)\,dr<\infty,
\end{eqnarray*}
where
$|K|$ is the volume of $K$ and $\sigma(\mathbb{S}^{N-1})$ is the area
of unit sphere.
\end{ex}

\section{Nakao integrals}\label{sec:Nakao}

Now we are in a position to define an extension of Nakao's
divergence-like CAF of zero energy and stochastic integrals with
respect to it in our setting.

Let $\mathcal{N}_c^*\subset\mathcal{N}_c$ denote the class of
continuous additive
functionals of the form $N^u+\int_0^\cdot g(X_s)\,ds$ for some $u\in
\mathcal{F}$ and $g\in L^2(E;m)$. Nakao \cite{Na} constructed a
linear operator
$\Gamma$ from $\accentset{\circ}{\mathcal{M}}$ into $\mathcal
{N}_c^*$ in the following
way: for every $Z\in
\accentset{\circ}{\mathcal{M}}$, there is a unique $w\in\mathcal
{F}$ such that
%
\begin{equation}\label{eq:nakaoriezse}
\mathcal{E}_1( w, f)= \tfrac{1}{2}\mu_{\langle M^f + M^{f, \kappa},
Z\rangle}(E)\qquad
\mbox{for every } f \in\mathcal{F}.
\end{equation}
This unique $w$ is denoted by $\gamma(Z)$.
The operator $\Gamma$ is defined by
%
\begin{equation}\label{eq:nakaoCAF}
\Gamma(Z)_t:=N_t^{\gamma(Z)} -\int_0^t \gamma(Z)(X_s) \,ds\qquad
\mbox{for } Z \in\accentset{\circ}{\mathcal{M}}.
\end{equation}
It is shown in Nakao \cite{Na} that $\Gamma(Z)$
can be characterized by the following equation:
%
\begin{equation}\label{eq:nakaoCAFries}
\lim_{t\downarrow0} \frac1{t}{\mathbb{E}}_{g \cdot m} [ \Gamma
(Z)_t ]
= - \frac12 \mu_{\langle M^g + M^{g, \kappa},   Z\rangle}(E)\qquad
\mbox{for every } g \in\mathcal{F}_b.
\end{equation}
Here $\mathcal{F}_b:=\mathcal{F}\cap L^\infty(E;m)$.
So, in particular, we have $\Gamma
(M^u)=N^u$ for $u\in\mathcal{F}$.
Moreover, we have the following:
\begin{lem}\label{lem:extendedNu} It holds that $\Gamma(M^u)=N^u$ for
$u\in\mathcal{F}_e$.
\end{lem}
\begin{pf}
Fix $u\in\mathcal{F}_e$ and let $\{u_n\}$ be an approximating
$\mathcal{E}$-Cauchy sequence
such that $u_n\to u$ $m$-a.e. In view of the proof of
Theorem 5.2.2 in \cite{FOT}, by taking a subsequence $\{n_k\}$,
$\{u_{n_k}(X_t)\}$, $M_t^{u_{n_k}}$ and $N_t^{u_{n_k}}$ uniformly
converges to
$u(X_t)$, $M_t^u$ and $N_t^u$, respectively,
on any compact subinterval of $[0,\infty[$ ${\mathbb{P}}_x$-a.s. for
q.e. $x\in E$.
From Theorem 3.2 in
\cite{Na}, by taking another subsequence, $\Gamma(M^{u_{n_k}})$
converges to
$\Gamma(M^u)$ uniformly on any finite interval ${\mathbb{P}}_x$-a.s. for
q.e. $x\in E$. Since $\Gamma(M^{u_{n_k}})=N^{u_{n_k}}$, we have
$\Gamma(M^u)=N^u$.
\end{pf}

In the same way of Nakao \cite{Na} (cf. (3.13) in \cite{CFKZStoch}),
we can define a stochastic integral
used by the operator $\Gamma$:
for $M\in\accentset{\circ}{\mathcal{M}}$ with its jump function
$\varphi
\in\accentset{\circ}{\mathcal{J}}$
and $f\in\mathcal{F}_e\cap L^2(E;\mu_{\langle M\rangle})$, we set
%
\begin{eqnarray}\qquad
\label{eq:nakaoint}
&&\int_0^t f (X_s)\,d\Gamma(M)_s\nonumber\\[-8pt]\\[-8pt]
&&\qquad:= \Gamma( f * M)_t
-\frac12
\langle M^{f, c}+M^{f, j},   M^{c}+M^{j}+K\rangle_t ,\qquad  t\in
[0,\infty[, \nonumber
\end{eqnarray}
where $(f* M)_t= \int_0^t f(X_{s-}) \,dM_s$ and $K\in\accentset{\circ}{\mathcal{M}}^d$ with
$K_t-K_{t-}=-\mathbf{1}_{E\times E}(\varphi+\overline{\varphi
})(X_{t-},X_t)$ $t\in\ ]0,\infty[$ ${\mathbb{P}}_x$-a.s. for
q.e. $x\in E$. Equation
(\ref{eq:nakaoint}) is well defined under ${\mathbb{P}}_x$ for
q.e. $x\in E$.
In this paper, we call the operator $\Gamma$ \textit{Nakao operator} and
the integral
(\ref{eq:nakaoint}) \textit{Nakao integral}.
\begin{rem}\label{rem:NakaoIntegral}
Equation (\ref{eq:nakaoint}) is consistent with the extension of
Nakao integral developed in
\cite{CFKZStoch} up to $\zeta$ under ${\mathbb{P}}_m$ (see
Theorem 6.3 in
\cite{KwstochII}).
\end{rem}

For any $M\in\accentset{\circ}{\mathcal{M}}^c_{
\mathrm{loc}}=\mathcal{M}_{\mathrm{loc}}^{c,\auf 0,\zeta\auf }$
(in particular for $M\in\accentset{\circ}{\mathcal{M}}
^c_{f\mbox{-}\mathrm{loc}}=\mathcal{M}_{\mathrm{loc}}^c$),
$\Gamma(M)$ can be defined as an element in ${\mathcal{N}}_{c,
f\mbox{-}\mathrm{loc}}$.
To see this, we need the following lemma extending Lemma 3.4 in \cite{Na}:
\begin{lem}[(Local property of $\Gamma$ on
$\accentset{\circ}{\mathcal{M}}
^c$)]\label{lem:divlikeCAFlocal}
Let $M\in\accentset{\circ}{\mathcal{M}}^c$ and $G$ be a (q.e.)
finely open
Borel set.
Suppose that $M_t=0$ ${\mathbb{P}}_m$-a.e. on $\{t<\tau_G\}$ for each $t>0$.
Then $\Gamma(M)_t=0$ for all $t\in[0,\sigma_{E\setminus G}]$
${\mathbb{P}}
_x$-a.s. for q.e. $x\in E$.
\end{lem}
\begin{pf}
By assumption, ${\mathbb{E}}_m[\langle M\rangle_{t\wedge\sigma
_{E\setminus G}}]=
{\mathbb{E}}_m[\langle M\rangle_{t\wedge\tau_{G}}]=
{\mathbb{E}}_m[M_{t\wedge\tau_{G}}^2]=0$ for each fixed $t>0$.
Then $\mu_{\langle M\rangle}(G)=0$ by Lemma 5.1.10(iii) in \cite{FOT}.
Let $h\in\mathcal{F}_F$ for a closed set $F$ with $F\subset G$.
Then $\mu_{\langle M^{h,c}\rangle}(E\setminus G)=0$ by
Lemma \ref{lem:MAFcontpartvanish}(1).
From this, we have
\[
\mu_{\langle M^{h,c},M\rangle}(E)^2\leq2\mu_{\langle M^{h,c}\rangle
}(E)\mu_{\langle M\rangle}(G)+2
\mu_{\langle M\rangle}(E)\mu_{\langle M^{h,c}\rangle}(E\setminus G)=0.
\]
Hence $\mathcal{E}_1(\gamma(M),h)=0$ for any $h\in\mathcal{F}_F$
with $F\subset G$.
Since $(\mathcal{E}_G,\mathcal{F}_G)$ is a quasi-regular Dirichlet
form on $L^2(G;m)$,
there exists an $\mathcal{E}_G$-nest $\{F_n\}$ of compact sets of $G$ (see
Lemma 3.4 in \cite{Kwfunc}).
From this, for any $h\in\mathcal{F}_G$, there exists
$h_k\in\bigcup_{n=1}^{\infty}\mathcal{F}_{F_n}$ such that $\{h_k\}$
$\mathcal{E}
^{1/2}_1$-converges to $h$ as
$k\to\infty$. Therefore, $\mathcal{E}_1(\gamma(M),h)=0$ for any
$h\in\mathcal{F}
_G$, which implies
$\Gamma(M)_t=N^{\gamma(M)}_t-\int_0^t\gamma(M)(X_s)\,ds=0$ for
$t<\sigma_{E\setminus G}$
by way of Lemma 5.4.2(ii) in \cite{FOT}.
\end{pf}

Let $({\mathcal{M}^{d}_{\mathrm{loc}}})_*$ [resp., $(\mathcal
{M}^{d}_{\mathrm{loc}})_{\mathrm{as}}$]
be the subclass of ${\mathcal{M}_{\mathrm{loc}}^d}$ associated to
$\mathcal{J}_*/\sim$ (resp., $\mathcal{J}_{\mathrm{as}}/\sim$) and
$({\mathcal{M}^{d,\auf 0,\zeta\auf }_{\mathrm{loc}}})_*$ [resp., $(\mathcal
{M}^{d,\auf 0,\zeta\auf }_{\mathrm{loc}})_{\mathrm{as}}$]
the subclass of $\mathcal{M}^{d,\auf 0,\zeta\auf }_{\mathrm{loc}}$ associated to
$\widehat{\mathcal{J}}_*/\sim$ (resp.,  $\widehat{\mathcal
{J}}_{\mathrm{as}}/\sim$).

We say that $M\in({\mathcal{M}_{\mathrm{loc}}})_*$ [resp., $M\in
({\mathcal{M}
_{\mathrm{loc}}})_{\mathrm{as}}$] if and only if its purely
discontinuous part
$M^d$ is in $({\mathcal{M}_{\mathrm{loc}}^d})_*$ [resp., $({\mathcal
{M}_{\mathrm{loc}}^d})_{\mathrm{as}}$], and
the classes $(\mathcal{M}^{\auf 0,\zeta\auf }_{\mathrm{loc}})_*$ and $(\mathcal
{M}^{\auf 0,\zeta
\auf }_{\mathrm{loc}})_{\mathrm{as}}$ are similarly defined.
For $M\in(\mathcal{M}^{\auf 0,\zeta\auf }_{\mathrm{loc}})_*$ with its jump function
$\varphi\in\widehat{\mathcal{J}}_*$,
let $M^c\in\mathcal{M}_{\mathrm{loc}}^{c,\auf 0,\zeta\auf }$ be its continuous
part and
take $M^j\in\mathcal{M}_{\mathrm{loc}}^j$ associated with $\mathbf
{1}_{E\times
E}\varphi\in\mathcal{J}_*$ and
$K\in\mathcal{M}^j_{\mathrm{loc}}$ constructed in Corollary \ref{cor:MAFreversion*}
associated with
$-\mathbf{1}_{E\times E}(\varphi+\overline{\varphi})\in\mathcal{J}_*$.


We shall extend $\Gamma$ over $(\mathcal{M}^{\auf 0,\zeta\auf }_{\mathrm{loc}})_*$ and
establish (\ref{eq:nakaoint}) for more general integrands and
integrators under
${\mathbb{P}}_x$ for q.e. $x\in E$. To do this we need the following lemma:
\begin{lem}\label{lem:localpropnkoint}
Take $M\in(\mathcal{M}^{\auf 0,\zeta\auf }_{\mathrm{loc}})_*$ with its jump function
$\varphi\in\widehat{\mathcal{J}}_*$.
Let $G$ be a (q.e.) finely open
Borel set satisfying $\mathbf{1}_{G\times E}\varphi,\mathbf
{1}_{G\times E}\overline{\varphi}\in\accentset{\circ}{\mathcal{J}}$.
Take $f\in\mathcal{F}_e\cap L^2(E;\mu_{\langle M\rangle})$ and
suppose that
$f=0$ $m$-a.e. on $G$. Then we have
$\Gamma(f*M)_t= \frac12\langle M^{f,c}+M^{f,j}, M^c+M^j+K\rangle_t$
for $t\in[0,\sigma_{E\setminus G}[$
${\mathbb{P}}_x$-a.s. for q.e. $x\in E$.
\end{lem}
\begin{pf} We show that for any $g\in(\mathcal{F}_G)_b$
%
\begin{equation} \label{eq:bddvarcaf}
\lim_{t\to0}\frac{1}{t}{\mathbb{E}}_{gm} \biggl[\Gamma(f*M)_t- \frac
12\langle
M^{f,c}+M^{f,j}, M^c+M^j+K\rangle_t \biggr]=0.
\end{equation}
Then we can obtain the assertion from (\ref{eq:bddvarcaf})
in view of the proof of Theorem 2.2 in \cite{Na} and Lemma 5.4.4 in
\cite{FOT}.
We know
\[
\lim_{t\to0}\frac{1}{t}{\mathbb{E}}_{gm} [\Gamma(f*M)_t ]=-\frac
{1}{2}\int
_Ef(x)\mu_{\langle M^g+M^{g,\kappa}, M\rangle}(dx).
\]
So for (\ref{eq:bddvarcaf}) it suffices to show
%
\begin{equation}\label{eq:bddvarder}\qquad
\int_Ef(x)\mu_{\langle M^g+M^{g,\kappa}, M\rangle}(dx)=-\int
_Eg(x)\mu_{\langle M^{f,c}+M^{f,j}, M^c+M^j+K\rangle}(dx).
\end{equation}
Noting $fg=0$ q.e. on $E$ and the derivation properties of continuous\break
part and jumping part of energy measures (see the proof of Lemma 3.1\break in
\cite{Na}), we see $\int_Ef \,d\mu_{\langle M^{g,c},  M^c\rangle
}+\int_Eg\,d\mu_{\langle M^{f,c},  M^c\rangle}=0$,
$\int_Ef \,d\mu_{\langle M^{g,j},  M^j\rangle}+\break\int_Eg\,d\mu_{\langle
M^{f,j},  M^j+K\rangle}=0$ and
$\int_Ef \,d\mu_{\langle M^{g,\kappa},  M^{\kappa}\rangle}=0$,
which imply (\ref{eq:bddvarder}).
\end{pf}
\begin{df}[(Extensions of Nakao operators and Nakao integrals)]
\label{df:localizationGamma}
Fix $M\in(\mathcal{M}^{\auf 0,\zeta\auf }_{\mathrm{loc}})_*$ with its jump function
$\varphi\in\widehat{\mathcal{J}}_*$
and $f\in\dot{\mathcal{F}}_{\mathrm{loc}}$.
Let $\{G_k\}\in\Theta$ be a common nest such that $\mu_{\langle
M\rangle}(G_k)<\infty$, $f=f_k$ $m$-a.e. on $G_k$ for some
$f_k\in\mathcal{F}_b$,
$\mathbf{1}_{G_k\times E}\varphi,\mathbf{1}_{G_k\times E}\overline
{\varphi}\in\accentset{\circ}{\mathcal{J}}$ for each $k\in\mathbb{N}$.
Set $E_k:=\{x\in E\mid{\mathbb{E}}_x[\int_0^{\tau
_{G_k}}e^{-t}g(X_t)\,dt]>1/k\}
$ for $g\in L^2(E;m)$ with $0<g\leq1$ $m$-a.e.
Then $e_k:=k{\mathbb{E}}_x[\int_0^{\tau_{G_k}}e^{-t}g(X_t)\,dt]\wedge
1\in\mathcal{F}_{G_k}$
satisfies $\mathbf{1}_{E_k}\leq e_{k}\leq\mathbf{1}_{G_k}$ q.e. on $E$.
In view of Lemma 3.3 in \cite{Kwfunc}, we have $\{E_k\}\in\Theta$.
\end{df}

We now define
\[
\Gamma(M)_t:=\Gamma(e_k*M)_t-\tfrac12\langle M^{e_k,c}+M^{e_k,j},
M^c+M^j+K\rangle_t \qquad\mbox{for }t\in[0,\sigma_{E\setminus E_k}[
\]
for each $k\in\mathbb{N}$
${\mathbb{P}}_x$-a.s. for q.e. $x\in E$. For $M\in(\mathcal{M}^
{\auf 0,\zeta\auf }_{\mathrm{loc}})_*$ and $f\in\mathcal{F}_e\cap L^2(E;\break\mu_{\langle M\rangle
})$, we set
\[
\int_0^tf(X_s)\,d\Gamma(M)_s:=\Gamma(f*M)_t-\frac12\langle
M^{f,c}+M^{f,j}, M^c+M^j+K\rangle_t
\]
for $t\in[0,\infty[$
${\mathbb{P}}_x$-a.s. for q.e. $x\in E$. For general $f\in\dot
{\mathcal{F}}_{\mathrm{loc}}$ and
$M\in(\mathcal{M}^{\auf 0,\zeta\auf }_{\mathrm{loc}})_*$ as above, we set
\[
\int_0^tf(X_s)\,d\Gamma(M)_s:=\int_0^t(fe_k)(X_s)\,d\Gamma(M)_s \qquad\mbox{for }t\in[0,
\sigma_{E\setminus E_k}[
\]
for each $k\in\mathbb{N}$
${\mathbb{P}}_x$-a.s. for q.e. $x\in E$. Note that $fe_k\in\mathcal
{F}_b\cap L^2(E;\mu
_{\langle M\rangle})$ for each $k\in\mathbb{N}$. These are well
defined for all $t\in[0,\infty[$ ${\mathbb{P}}_x$-a.s. for
q.e. $x\in E$
in view of Lemma \ref{lem:localpropnkoint} and are elements in
${\mathcal{N}}_{c,f\mbox{-}\mathrm{loc}}$.

For $f\in\dot{\mathcal{F}}_{\mathrm{loc}}$ and $M\in(\mathcal{M}^{\auf 0,\zeta\auf }_{\mathrm{loc}})_*$, we see
%
\begin{equation}\label{eq:contnakaoint}
\int_0^tf(X_s)\,d\Gamma(M^c)_s=\Gamma(f*M^c)_t-\frac{1}{2}\langle
M^{f,c},M^c\rangle_t
\end{equation}
for all $t\in[0,\infty[$ ${\mathbb{P}}_x$-a.s. for q.e. $x\in E$, where
$\Gamma(f*M^c)_t$ can be defined
by way of Lemma \ref{lem:divlikeCAFlocal}.
\begin{rem}
(1)
In \cite{CFKZStoch}, we define extensions of Nakao operators and
Nakao

{\smallskipamount=0pt
\begin{enumerate}[(2)]
\item[ ]
integrals in terms of time reverse
operators, which are defined up to $\zeta$ under~${\mathbb{P}}_m$,
and the Nakao integral is also refined for integrator $\Gamma(M)$,
$M\in\accentset{\circ}{\mathcal{M}}$ and integrand $f(X)$ for
$f\in\dot
{\mathcal{F}}_{\mathrm{loc}}$ under ${\mathbb{P}}_x$ for q.e. $x\in E$.
So the Nakao integral in
Definition \ref{df:localizationGamma} is a pure extension of this refinement.
Though the condition on the integrator of our Nakao
integrals is rather restrictive than theirs described to be up to
$\zeta$ under ${\mathbb{P}}_m$,
it is defined for all time under the law
for quasi-everywhere starting points.
\item[(2)] The extensions of Nakao operators and Nakao
integrals in \cite{CFKZStoch} are consistent with
our corresponding notions up to $\zeta$ under ${\mathbb{P}}_m$ (see
Theorem 6.3 in \cite{KwstochII}).
\end{enumerate}}
\end{rem}

The following lemma is needed to establish the generalized It\^o formula.
\begin{lem}[(Local property of extended Nakao integral)]\label{lem:localNakaoInt}
Take
\[
M\in\bigl(\mathcal{M}^{\auf 0,\zeta\auf }_{\mathrm{loc}}\bigr)_*
\]
with its jump function
$\varphi\in\widehat{\mathcal{J}}_*$
and $f\in\dot{\mathcal{F}}_{\mathrm{loc}}$. Let $G$ be a (q.e.) finely open
Borel set.
Suppose that $f=0$ $m$-a.e. on $G$. Then
\[
\int_0^t f(X_s)\,d\Gamma(M)_s=0
\]
holds for all $t\in[0,\sigma_{E\setminus G}[$ ${\mathbb{P}}_x$-a.s. for
q.e. $x\in E$.
\end{lem}
\begin{pf}
Let $\{E_k\}\in\Theta$ and $e_k\in\mathcal{F}$ be constructed as in
Definition \ref{df:localizationGamma}.
Since $fe_k\in\mathcal{F}_b\cap L^2(E;\mu_{\langle M\rangle})$ and
$\mathbf
{1}_{(G\cap E_k)\times E}\varphi,
\mathbf{1}_{(G\cap E_k)\times E}\overline{\varphi}\in\accentset{\circ}{\mathcal{J}}$,
we can apply Lem\-ma~\ref{lem:localpropnkoint}
so that
\[
\int_0^t (fe_k)(X_s)\,d\Gamma(M)_s=0  \qquad\mbox{holds for } t\in
[0,\sigma_{E\setminus G}\wedge\sigma_{E\setminus E_k}[,
\]
${\mathbb{P}}_x$-a.s. for q.e. $x\in E$. Therefore we obtain the
desired assertion.
\end{pf}

The following propositions are an addendum (cf. Theorems 4.1 and 4.2 in
\cite{CFKZStoch}).
We omit its proofs.
\begin{pr}
Take $M\in(\mathcal{M}^{\auf 0,\zeta\auf }_{\mathrm{loc}})_*$ and $f\in\dot
{\mathcal{F}}_{\mathrm{loc}}$.
Suppose that $\Gamma(M)$ is a CAF $A$ of finite variation on $\auf 0,\zeta\auf $.
Then
\[
\int_0^tf(X_s)\,d\Gamma(M)_s=\int_0^tf(X_s)\,dA_s
\]
holds for all $t\in[0,\zeta[$ (for all $t\in[0,\infty[$ provided
$M\in(\mathcal{M}_{\mathrm{loc}})_*$) ${\mathbb{P}}_x$-a.s. for q.e. $x\in E$.
\end{pr}
\begin{pr}
Take $M\in(\mathcal{M}^{\auf 0,\zeta\auf }_{\mathrm{loc}})_*$ and $f,g\in\dot
{\mathcal{F}}_{\mathrm{loc}}$.
Then
\[
\int_0^tg(X_s)\,d \biggl(\int_0^{\cdot}f(X_u)\,d\Gamma(M)_u \biggr)_s=\int
_0^t(fg)(X_s)\,d\Gamma(M)_s
\]
holds for all $t\in[0,\zeta[$ (for all $t\in[0,\infty[$ provided
$M\in(\mathcal{M}_{\mathrm{loc}})_*$) ${\mathbb{P}}_x$-a.s. for q.e. $x\in E$.
\end{pr}

\section{Stochastic integrals by Dirichlet processes}\label{sec:StochasticIntegral}

\begin{df}[(Dirichlet processes)]
For $M\in(\mathcal{M}^{\auf 0,\zeta\auf }_{\mathrm{loc}})_*$,
we set $A:=M+\Gamma(M)$, which can be defined by way of
Definition \ref{df:localizationGamma}.
Note that $A$ is defined on $[0,\infty[$ ${\mathbb{P}}_x$-a.s. for
q.e. $x\in E$
if $M\in(\mathcal{M}_{\mathrm{loc}})_*$; otherwise,
it is only defined on $[0,\zeta[$ ${\mathbb{P}}_x$-a.s. for
q.e. $x\in E$.
For $M=M^u$ with $u\in\mathcal{F}_e$, we see $A=A^u$, where
$A_t^u:=u(X_t)-u(X_0)$. For $M\in\accentset{\circ}{\mathcal{M}}$
and each
$i=c,d,j,\kappa$, we further set
$A^{i}_t:=M^{i}_t+\Gamma(M^{i})_t$
and write $A^{u,i}:=A^i$ if $M=M^u$, $u\in\mathcal{F}_e$.
\end{df}

We see $A=A^c+A^d=A^c+A^j+A^{\kappa}$ for $A=M+\Gamma(M)$, $M\in
\accentset{\circ}{\mathcal{M}}$.
By (\ref{eqn:fuku}) and Lemma \ref{lem:extendedNu},
we have $A^u=A^{u,c}+A^{u,d}=A^{u,c}+A^{u,j}+A^{u,\kappa}$ for $u\in
\mathcal{F}_e$.
\begin{df}[(Stochastic integrals by Dirichlet processes)]
\label{df:FiskStratonovich}
Take and fix $M\in\accentset{\circ}{\mathcal{M}}$. For $f\in
L^2(E;\mu
_{\langle M^c\rangle})$, we set
\[
\int_0^t f(X_s)\circ dA_s^c:=(f*M^c)_t+\Gamma(f*M^c)_t
\]
for all $t\in[0,\infty[$ ${\mathbb{P}}_x$-a.s. for q.e. $x\in E$.
For
$f\in\mathcal{F}_e\cap L^2(E;\mu_{\langle M\rangle})$, we set
\begin{eqnarray*}
\int_0^t f(X_s)\circ dM_s&:=&(f*M)_t+\frac{1}{2}[M^{f},M]_t, \\
\int_0^t f(X_{s-})  \,dA_s&:=& (f*M)_t+\int_0^t f(X_s)\,d\Gamma(M)_s,\\
\int_0^t f(X_s)\circ dA_s&:=& \int_0^t f(X_s)\circ dM_s+\int_0^t
f(X_s)\,d\Gamma(M)_s
\end{eqnarray*}
for all $t\in[0,\infty[$ ${\mathbb{P}}_x$-a.s. for q.e. $x\in E$.
Recall $(f*M)_t=\int_0^tf(X_{s-})\,dM_s$ for $f\in\mathcal{F}_e\cap
L^2(E;\mu
_{\langle M\rangle})$
(see the proof of Lemma \ref{lem:Riemannsumapproxmar}).
We call
$(f*M)_t$ [resp., $\int_0^tf(X_{s-})  \,dA_s$]
the \textit{It\^o integral} of $f(X)$ with integrator $M$ (resp., $A$)
and
$\int_0^tf(X_s)\circ dM_s$ [resp., $\int_0^tf(X_s)\circ dA_s$]
the \textit{Fisk--Stratonovich integral} of $f(X)$ with integrator $M$
(resp., $A$).
\end{df}
\begin{rem}\label{rem:FiskStratonovich}
%
(1)
For the definition of $\int_0^tf(X_s)\circ dM_s^c$ for
$f\in\mathcal{F}
_e\cap L^2(E;$%

{\smallskipamount=0pt
\begin{enumerate}[(2)]
\item[ ] $\mu_{\langle M^c\rangle})$,
we need $f\in\mathcal{F}_e$, which is unnecessary to define $\int
_0^tf(X_s)\circ dA_s^c$.
\item[(2)] Our definitions of Fisk--Stratonovich-type integrals are somewhat
different from what is
found in Protter \cite{ProtStoch} or in Chapter VI of Meyer \cite
{Meyerstochs} except for the case of diffusions.
The definition of $\int_0^t f(X_s)\circ dM_s$ in \cite{ProtStoch} or
\cite{Meyerstochs}
is given by $(f*M)_t+
\frac12\langle M^{f,c},M^{c}\rangle_t$, which has an advantage
to give a version of It\^o's formula in terms of their
Fisk--Stratonovich integrals
(see II. Theorem~34, V. Theorems 20 and 21 in \cite{ProtStoch}),
but it only has a Riemann-sum approximation
under that $f(X_s)$ and $u(X_s)$ have no jumps in common (see V.
Theorem 26 in~\cite{ProtStoch}).
Our definition of $\int_0^t f(X_s)\circ dM_s$
admits such an approximation in the framework of semi-martingales at
least (cf.
Definition 3.9.21 in \cite{BichtelerStoch} or Problems 9.12 and 9.13
in \cite{HWY}).
On the other hand, Kurtz, Pardoux and Protter \cite{KPP}
give a different definition for Fisk--Stratonovich-type integrals provided
the underlying process is a solution of an SDE driven by semimartingales.
Our definitions are also different from theirs.
\end{enumerate}}
\end{rem}

Now take a jump function $\varphi\in{\mathcal{J}}$ associated to
a given $M\in{\mathcal{M}}_{\mathrm{loc}}$. We set for each $\ell\in\mathbb{N}$
\[
M_t^{d,\ell}:=\sum_{0<s\leq t}\varphi_{\ell}(X_{s-},X_s) -
\int_0^tN(\varphi_{\ell})(X_s)\,dH_s.
\]
In the same way of the proof of Lemma
\ref{lem:purediscontimartincorrespondince}, if
$M\in\accentset{\circ}{\mathcal{M}}$, then
\[
\mathbf{e}(M^{d}-M^{d,\ell})=\frac12
\int_E\int_{E_{\partial}}\varphi^2(x,y)\mathbf{1}_{\{|\varphi
(x,y)|\leq1/\ell\}}N(x,dy)\mu_H(dx).
\]
The stochastic integrals $f*M^d$ and $f*M^{d,\ell}$ for $M\in
{\mathcal{M}
}_{\mathrm{loc}}$ with
$f\in\dot{\mathcal{F}}_{\mathrm{loc}}\cap L^2(E;\mu_{\langle M^{d}\rangle
})$ and
$f(\partial)=0$ belong to $\accentset{\circ}{\mathcal{M}}$, and
satisfy that
\[
(f*M^{d,\ell})_t=\sum_{0<s\leq t}f(X_{s-})\varphi_{\ell}(X_{s-},X_s)
-\int_0^tf(X_s)N(\varphi_{\ell})(X_s)\,dH_s
\]
holds for all $t\in[0,\infty[$ ${\mathbb{P}}_x$-a.s. for q.e. $x\in
E$ and
\begin{eqnarray*}
&&
\mathbf{e}(f*M^{d}-f*M^{d,\ell})\\
&&\qquad=\frac12
\int_Ef^2(x)\int_{E_{\partial}}\varphi^2(x,y)\mathbf{1}_{\{
|\varphi(x,y)|\leq1/\ell\}}N(x,dy)\mu_H(dx).
\end{eqnarray*}
\begin{lem}\label{lem:compensator2}
(1)
Take $M\in{\mathcal{M}}_{\mathrm{loc}}$ with its jump function
$\varphi\in{\mathcal{J}}$. Then for

{\smallskipamount=0pt
\begin{enumerate}[(2)]
\item[ ]
$g\in\dot{\mathcal
{F}}_{\mathrm{loc}}\cap L^2(E;\mu
_{\langle M^d\rangle})$ with
$g(\partial)=0$,
%
\begin{eqnarray}\label{eq:compensator1}
{\Gamma}(g*M^{d,\ell})_t&=&\frac{1}{2}\int_0^tN\bigl(\mathbf{1}_{E\times
E}(g\varphi_{\ell}-\overline{g\varphi_{\ell}})\bigr)(X_s)\,dH_s
\nonumber\\[-8pt]\\[-8pt]
&&{}+\int_0^tg(X_s)\varphi_{\ell}(X_s,\partial)N(X_s,\{
\partial\})\,dH_s
\nonumber
\end{eqnarray}
holds for all $t\in[0,\infty[$ ${\mathbb{P}}_x$-a.s. for q.e. $x\in E$.
Moreover, for
$f\in\mathcal{F}_e\cap L^2(E;\mu_{\langle g*M^{d,\ell}\rangle})$,
%
\begin{eqnarray}\label{eq:compensator2}
&&\int_0^tf(X_s)\,d{\Gamma}(g*M^{d,\ell})_t\nonumber\\
&&\qquad=\frac12
\int_0^tf(X_s)N\bigl(\mathbf{1}_{E\times E}(g\varphi_{\ell}-\overline
{g\varphi_{\ell}})\bigr)(X_s)\,dH_s\\
&&\qquad\quad{}+
\int_0^tf(X_s)g(X_s)\varphi_{\ell}(X_s,\partial)N(X_s,\{\partial\}
)\,dH_s\nonumber
\end{eqnarray}
holds for all $t\in[0,\infty[$ ${\mathbb{P}}_x$-a.s. for q.e. $x\in E$.
More generally if $M\in\mathcal{M}^{\auf 0,\zeta\auf }_{\mathrm{loc}}$ with its
jump function
$\varphi\in\widehat{\mathcal{J}}$, then for $M^{j,\ell}\in
\mathcal{M}_{\mathrm{loc}}^j$ with its
jump function $\mathbf{1}_{E\times E}\varphi_{\ell}\in\mathcal{J}$ and
$g\in\dot{\mathcal{F}}_{\mathrm{loc}}\cap L^2(E;\mu_{\langle M^{j,\ell
}\rangle})$
with $g(\partial)=0$,
%
\begin{equation}\label{eq:compensator1*}
{\Gamma}(g*M^{j,\ell})_t=\frac{1}{2}\int_0^tN\bigl(\mathbf{1}_{E\times
E}(g\varphi_{\ell}-\overline{g\varphi_{\ell}})\bigr)(X_s)\,dH_s
\end{equation}
holds for all $t\in[0,\infty[$ ${\mathbb{P}}_x$-a.s. for q.e. $x\in
E$, and for
$f\in\mathcal{F}_e\cap L^2(E;\break\mu_{\langle g*M^{j,\ell}\rangle})$,
%
\begin{equation}\label{eq:compensator2*}
\int_0^tf(X_s)\,d{\Gamma}(g*M^{j,\ell})_t=\frac12
\int_0^tf(X_s)N\bigl(\mathbf{1}_{E\times E}(g\varphi_{\ell}-\overline
{g\varphi_{\ell}})\bigr)(X_s)\,dH_s
\end{equation}
holds for all $t\in[0,\infty[$ ${\mathbb{P}}_x$-a.s. for q.e. $x\in E$.
\item[(2)]
Take $M\in(\mathcal{M}^{\auf 0,\zeta\auf }_{\mathrm{loc}})_*$ and $K^{\ell}\in
\mathcal{M}
_{\mathrm{loc}}^j$ associated with
$-\mathbf{1}_{E\times E}(\varphi_{\ell}+\overline{\varphi_{\ell
}})\in
\mathcal{J}_*$.
Then for $g\in\dot{\mathcal{F}}_{\mathrm{loc}}\cap L^2(E;\mu_{\langle
K^{\ell
}\rangle})$ with $g(\partial)=0$, we have
%
\begin{equation}\label{eq:compensator3}
{\Gamma}(g*K^{\ell})_t=\frac{1}{2}\int_0^tN\bigl(\mathbf{1}_{E\times
E}(\bar{g}-g)(\varphi_{\ell}
+\overline{\varphi_{\ell}})\bigr)(X_s)\,dH_s
\end{equation}
holds for all $t\in[0,\infty[$ ${\mathbb{P}}_x$-a.s. for q.e. $x\in E$.
Moreover, for
$f\in\mathcal{F}_e\cap L^2(E;\mu_{\langle g*K^{\ell}\rangle})$,
%
\begin{eqnarray}\label{eq:compensator4}
&&\int_0^tf(X_s)\,d{\Gamma}(g*K^{\ell})_t\nonumber\\[-8pt]\\[-8pt]
&&\qquad=\frac12
\int_0^tf(X_s)N\bigl(\mathbf{1}_{E\times E}(\bar{g}-g)(\varphi_{\ell}
+\overline{\varphi_{\ell}})\bigr)(X_s)\,dH_s\nonumber
\end{eqnarray}
holds for all $t\in[0,\infty[$ ${\mathbb{P}}_x$-a.s. for q.e. $x\in E$.
\end{enumerate}}
\end{lem}
\begin{cor}\label{cor:compasator2*}
Take $M\in(\mathcal{M}^{\auf 0,\zeta\auf }_{\mathrm{loc}})_*$ and $f\in\dot
{\mathcal{F}
}_{\mathrm{loc}}$. Let $K$ be an element in
$(\mathcal{M}^d_{\mathrm{loc}})_*$ associated with $-\mathbf{1}_{E\times
E}(\varphi
+\overline{\varphi})$ constructed in
Corollary \ref{cor:MAFreversion*}. Then we have that
\[
\int_0^tf(X_s)\,d\Gamma(K)_s=0
\]
holds for all $t\in[0,\infty[$ ${\mathbb{P}}_x$-a.s. for q.e. $x\in
E$. In
particular, $\Gamma(K)_t=0$
for all $t\in[0,\infty[$ ${\mathbb{P}}_x$-a.s. for q.e. $x\in E$.
\end{cor}
\begin{cor}\label{cor:compensator3}
Take $M\in\accentset{\circ}{\mathcal{M}}$ with its jump function
$\varphi\in\accentset{\circ}{\mathcal{J}}$.
Let $K^{\ell}$ be an element in
$\accentset{\circ}{\mathcal{M}}^j$ associated with $-\mathbf
{1}_{E\times
E}(\varphi_{\ell}+\overline{\varphi_{\ell}})$.
Set $\overline{A}{}^{d,\ell}:=A^{d,\ell}+\frac{1}{2}K^{\ell}$.
Then we have
\[
\overline{A}{}^{d,\ell}_t=\frac{1}{2}\sum_{0<s\leq
t}(\varphi_{\ell}-\overline{\varphi_{\ell}})(X_{s-},X_s)\mathbf
{1}_{\{
s<\zeta\}}+\varphi_{\ell}(X_{\zeta-},\partial)\mathbf{1}_{\{t\geq
\zeta\}}
\]
holds for all $t\in[0,\infty[$ ${\mathbb{P}}_x$-a.s. for q.e. $x\in E$.
Moreover, for
$f\in\mathcal{F}_e\cap L^2(E;\break\mu_{\langle M^{d}\rangle})$,
\begin{eqnarray*}
\int_0^tf(X_{s-})\, d\overline{A}{}^{d,\ell}_s&=&
\sum_{0<s\leq t}f(X_{s-})
\frac{\varphi_{\ell}-\overline{\varphi_{\ell
}}}{2}(X_{s-},X_s)\mathbf
{1}_{\{s<\zeta\}} \\
&&{}
+
f(X_{\zeta-})\varphi_{\ell}(X_{\zeta-},\partial)\mathbf{1}_{\{
t\geq\zeta\}}
\end{eqnarray*}
and
\begin{eqnarray*}
\int_0^t  f(X_{s-})\circ d\overline{A}{}^{d,\ell}_s&=&\sum_{0<s\leq
t}\frac{f(X_s) + f(X_{s-})}{2}
\frac{\varphi_{\ell}-\overline{\varphi_{\ell
}}}{2}(X_{s-},X_s)\mathbf
{1}_{\{s<\zeta\}} \\
&&{}
+f(X_{\zeta-})\varphi_{\ell}(X_{\zeta-},\partial)\mathbf{1}_{\{
t\geq\zeta\}}
\end{eqnarray*}
hold for all $t\in[0,\infty[$ ${\mathbb{P}}_x$-a.s. for q.e. $x\in E$.
\end{cor}
\begin{pf*}{Proof of Lemma \ref{lem:compensator2}}
We only prove (\ref{eq:compensator1}).
The proofs of (\ref{eq:compensator2}), (\ref{eq:compensator3}) and
(\ref{eq:compensator4}) are similar. Equation (\ref{eq:compensator1*})
[resp., (\ref{eq:compensator2*})] is clear from (\ref
{eq:compensator1}) [resp., (\ref{eq:compensator2})].
By~(\ref{eq:nakaoCAFries}), for $h \in\mathcal{F}_b$
\begin{eqnarray*}
&&\lim_{t\downarrow0} \frac1{t}{\mathbb{E}}_{h \cdot m} [ {\Gamma}
(g*M^{d,\ell})_t ]\\
&&\qquad= - \frac12\int_Eg  \,d\mu_{\langle M^h + M^{h, \kappa},
M^{d,\ell}\rangle}\\
&&\qquad=-\frac12\int_E g \,d\mu_{\langle M^{h,j}, M^{j,\ell}\rangle}
-\int_E g  \,d
\mu_{\langle M^{h,\kappa}, M^{\kappa,\ell}\rangle}
\\
&&\qquad=\int_{E\times E}h(x)(g\varphi_{\ell}-\overline
{g\varphi_{\ell}})(x,y)J(dx\,dy)
+\int_Eh(x)g(x)\varphi_{\ell}(x,\partial)\kappa(dx)\\
&&\qquad=\frac{1}{2}\int_E h N\bigl(\mathbf{1}_{E\times E}(g\varphi_{\ell
}-\overline{g\varphi_{\ell}})\bigr)\,d\mu_H
+\int_Ehg N\bigl(\mathbf{1}_{E\times\{\partial\}}\varphi_{\ell}\bigr)\,d\mu_H.
\end{eqnarray*}
Therefore, by Theorem 2.2 in \cite{Na}, we have the desired assertion.
\end{pf*}
\begin{pf*}{Proof of Corollary \ref{cor:compasator2*}}
Let $\{G_k\}\in\Theta$ be a common nest such that $f|_{G_k}\in
\mathcal{F}
_b|_{G_k}$ and
$\mathbf{1}_{G_k\times E}\varphi$, $\mathbf{1}_{G_k\times
E}\overline
{\varphi}\in\accentset{\circ}{\mathcal{J}}$ for each $k\in\mathbb{N}$.
Let $\{E_k\}\in\Theta$ be the nest
and let $e_k$ be the function constructed through $\{G_k\}$ as in
Definition \ref{df:localizationGamma}. Replacing $f$ with $fe_k$,
it suffices to prove the assertion for the case $f\in(\mathcal
{F}_{G_k})_b$ in
view of
Definition \ref{df:localizationGamma}.
For $f\in(\mathcal{F}_{G_k})_b$, we have that
\[
\int_0^tf(X_s)\,d\Gamma(K)_s=\Gamma(f*K)_t+\frac{1}{2}\langle
M^{f,j},K\rangle_t
\]
holds for $t\in[0,\infty[$ ${\mathbb{P}}_x$-a.s. for q.e. $x\in E$.
From Lemma \ref{lem:compensator2}, we have
$\Gamma(f*K^{\ell})_t+\frac{1}{2}\langle M^{f,j},K^{\ell}\rangle
_t=0$ holds for $t\in[0,\infty[$ ${\mathbb{P}}_x$-a.s. for
q.e. $x\in E$ in
view of Theorem 2.2 in \cite{Na}.
On the other hand, we see $f*K,f*K^{\ell},\mathbf{1}_{G_j}*K,\mathbf
{1}_{G_j}*K^{\ell}\in\accentset{\circ}{\mathcal{M}}$ for
$j>k$ with
\[
\mathbf{e}\bigl(f*(K-K^{\ell})\bigr)\leq\|f\|_{\infty}^2 \mathbf{e}\bigl(\mathbf
{1}_{G_j}*(K-K^{\ell})\bigr)\to0, \qquad \ell\to\infty.
\]
Hence we obtain the assertion in view of Theorem 5.2.1 in \cite{FOT}
and Theorem 3.2 in \cite{Na}.
\end{pf*}
\begin{pf*}{Proof of Corollary \ref{cor:compensator3}}
Since $K^{\ell}\in\accentset{\circ}{\mathcal{M}}$, we have from
Corollary \ref{cor:compasator2*} that
$\Gamma(K^{\ell})_t=0$ holds for $t\in[0,\infty[$ ${\mathbb{P}}_x$-a.s. for
q.e. $x\in E$.
Note that $M^{d,\ell}+\frac{1}{2}K^{\ell}$ is given by
\begin{eqnarray*}
&&\sum_{0<s\leq t}\frac{\varphi_{\ell}-\overline{\varphi_{\ell
}}}{2}(X_{s-},X_s)\mathbf{1}_{\{s<\zeta\}}+\varphi_{\ell}(X_{\zeta
-},\partial)\mathbf{1}_{\{t\geq\zeta\}} \\
&&\qquad{} -\int_0^tN \biggl(\mathbf{1}_{E\times E}\frac{\varphi_{\ell}-\overline
{\varphi_{\ell}}}{2} \biggr)(X_s)\,dH_s-\int_0^t\varphi_{\ell
}(X_s,\partial)N(X_s,\{\partial\})\,dH_s.
\end{eqnarray*}
Then we obtain the assertion in view of Lemma \ref{lem:compensator2}.
\end{pf*}
\begin{df}[(Extensions of stochastic integrals by Dirichlet
processes)]\label{df:extenstochinteg}
For $M\in(\mathcal{M}_{\mathrm{loc}})_*$ with its jump function
$\varphi\in\mathcal{J}_*$ and $f\in L^2(G_n;$ $\mu_{\langle
M^c\rangle})$ for each $n\in\mathbb{N}$ and some $\{G_n\}\in\Theta
$, we define
\[
\int_0^tf(X_s)\circ dA^c_s:=(f*M^c)_t+\Gamma(f*M^c)_t
\]
for all $t\in[0,\infty[$ ${\mathbb{P}}_x$-a.s. for q.e. $x\in E$.
Moreover, if $f\in\dot{\mathcal{F}}_{\mathrm{loc}}$ with $f(\partial)=0$, we define
\begin{eqnarray*}
\int_0^t f(X_{s-})\circ dM_s&:=&\int_0^tf(X_{s-})\,dM_s+\frac
{1}{2}[f(X),M]_t, \\
\int_0^t f(X_{s-})\,dA_s&:=&\int_0^tf(X_{s-})\,dM_s+\int_0^tf(X_s)\,d\Gamma
(M)_s,\\
\int_0^t f(X_{s})\circ dA_s&:=&\int_0^t f(X_{s-})\circ dM_s+\int
_0^tf(X_s)\,d\Gamma(M)_s
\end{eqnarray*}
for all $t\in[0,\infty[$ ${\mathbb{P}}_x$-a.s. for q.e. $x\in E$,
where $[f(X),M]_t :=\langle M^{f,c},M^c\rangle_t+\sum_{0<s\leq
t}(f(X_s)-f(X_{s-}))(M_s-M_{s-})$.
For $M\in(\mathcal{M}^{\auf 0,\zeta\auf }_{\mathrm{loc}})_*$ with its jump function
$\varphi\in\widehat{\mathcal{J}}_*$ and $f$ as above,
these are defined for all $t\in[0,\zeta[$ ${\mathbb{P}}_x$-a.s. for
q.e. $x\in E$.
We can define $\int_0^t f(X_{s-})\,d\overline{A}_s$, $\int_0^t
f(X_{s})\circ d\overline{A}_s$ for
$\overline{A}:=A+\frac{1}{2}K$ by replacing $M$ with $M+\frac{1}{2}K$.
Note that $\Gamma(K)\equiv0$.
\end{df}

Hereafter we use the following convention: let $f\in\dot{\mathcal{F}}_{\mathrm{loc}}$
with $f(\partial)=0$
and take $\phi,\psi\dvtx E_{\partial}\times E_{\partial}\to\mathbb{R}$
vanishing on
the diagonal such
that $|\phi|\leq M|\psi|$ on $E\times E_{\partial}$ for some $M>0$ and
$\sum_{0<s\leq t}\psi^2(X_{s-},X_s)<\infty$ for all $t\in[0,\infty
[$ ${\mathbb{P}}_x$-a.s. for q.e. $x\in E$.
If there exists a nest $\{G_n\}\in\Theta$
such that $f|_{G_n}\in\mathcal{F}_b|_{G_n}$ for each $n\in\mathbb{N}$
and a subsequence $\{\ell_k\}$ depending only on
$\{G_n\}$, $f$, $\phi$ and $\psi$
such that
$t\mapsto\sum_{0<s\leq t}f(X_{s-})\phi(X_{s-},X_s)
\mathbf{1}_{\{|\psi(X_{s-},X_s)|>1/\ell_k\}}
$
converges uniformly on each compact subinterval of $[0,\sigma
_{E\setminus G_n}[$ for each $n\in\mathbb{N}$ as $k\to\infty$
${\mathbb{P}}_x$-a.s. for q.e. $x\in E$,
then
we shall denote its limit by
\[
{\sum_{0<s\leq t}}^{\!\!\!\!  *}f(X_{s-})\phi(X_{s-},X_s).
\]
Note that if $t\mapsto\sum_{s\leq t}f(X_{s-})\phi(X_{s-},X_s)$
absolutely converges
uniformly on each compact subinterval of $[0,\infty[$ ${\mathbb{P}}
_x$-a.s. for q.e. $x\in E$, then
we can eliminate the symbol $*$ from the above sum.
We shall use $\sum_{s\leq t}^*(f(X_s)+f(X_{s-}))\phi(X_{s-},X_s)$ and
$\sum_{s\leq t}^*\phi(X_{s-},X_s)$ in a similar fashion.

We then have the following:
\begin{theorem}\label{thm:Stratonovichlocalrepre}
Let $f\in\dot{\mathcal{F}}_{\mathrm{loc}}$ and suppose that $f$ is extended to
be a
real-valued function on $E_{\partial}$
with $f(\partial)=0$.
Take $M\in(\mathcal{M}_{\mathrm{loc}})_*$ with its jump function
$\varphi\in\mathcal{J}_*$
and set $\overline{A}:=A+\frac{1}{2}K=M+\Gamma(M)+\frac{1}{2}K$, where
$K\in(\mathcal{M}_{\mathrm{loc}})_*$ associated with $-\mathbf{1}_{E\times
E}(\varphi
+\overline{\varphi})\in\mathcal{J}_*$. Then
\begin{eqnarray*}
\overline{A}_t&=&A_t^c+{\sum_{0<s\leq t}}^{\!\!\!\!*}
\mathbf{1}_{E\times E}\frac{\varphi-\overline{\varphi
}}{2}(X_{s-},X_s)+\varphi(X_{\zeta-},\partial)\mathbf{1}_{\{t\geq
\zeta\}},\\
\int_0^tf(X_{s-})  \,d\overline{A}_s&=&\int_0^tf(X_{s-})  \,d{A}_s^c\\
&&{}+
{\sum_{0<s\leq t}}^{\!\!\!\!*}f(X_{s-})
\mathbf{1}_{E\times E}\frac{\varphi-\overline\varphi
}{2}(X_{s-},X_s)\\
&&{}+f(X_{\zeta-})\varphi(X_{\zeta-},\partial)\mathbf
{1}_{\{t\geq\zeta\}}
\end{eqnarray*}
and
\begin{eqnarray*}
\int_0^tf(X_{s})\circ d\overline{A}_s&=&\int_0^tf(X_{s})\circ
\,dA_s^c\\
&&{}+
{\sum_{0<s\leq t}}^{\!\!\!\!*}\frac{f(X_s)+f(X_{s-})}{2}
\mathbf{1}_{E\times E}\frac{\varphi-\overline{\varphi
}}{2}(X_{s-},X_s)\\
&&{}+f(X_{\zeta-})\varphi(X_{\zeta-},\partial)\mathbf
{1}_{\{t\geq\zeta\}}
\end{eqnarray*}
hold for all $t\in[0,\infty[$ ${\mathbb{P}}_x$-a.s. for q.e. $x\in E$.
More generally, for $M\in(\mathcal{M}^{\auf 0,\zeta\auf }_{\mathrm{loc}})_*$
with its
jump function
$\varphi\in\widehat{\mathcal{J}}_*$, these expressions hold for
$t\in[0,\zeta[$
${\mathbb{P}}_x$-a.s. for q.e. $x\in E$.
\end{theorem}
\begin{pf}
First we assume $M\in\accentset{\circ}{\mathcal{M}}$ and $f\in
\mathcal{F}_e\cap
L^2(E;\mu_{\langle M\rangle})$.
Since $f*M^{d,\ell}$ converges $f*M^{d}$ in $(\accentset{\circ}{\mathcal{M}},\mathbf{e})$
as $\ell\to\infty$, there exists a common subsequence $\{\ell_k\}$
such that
$f*M^{d,\ell_k}$ [resp., $\Gamma(f*M^{d,\ell_k})$] uniformly
converges to $f*M^{d}$ [resp., $\Gamma(f*M^{d})$] on each compact
subinterval of $[0,\infty[$ ${\mathbb{P}}_x$-a.s. for q.e. $x\in E$
by Theorem 5.2.1 in \cite{FOT} and Theorem 3.2 in \cite{Na}.
On the other hand, $M^{d,\ell}$ also converges $M^{d}$ in
$(\accentset{\circ}{\mathcal{M}},\mathbf{e})$,
which yields that\vspace*{1pt} there exists a subsequence $\{\ell_k\}$ such that
$[M^{f,d},M^{d,\ell_k}]$ (resp., $\langle M^{f,d},M^{d,\ell_k}\rangle
$) uniformly converges to $[M^{f,d},M^{d}]$
(resp., $\langle M^{f.d},M^{d}\rangle$) on each compact subinterval of
$[0,\infty[$ ${\mathbb{P}}_x$-a.s. for q.e. $x\in E$.
Therefore, for such subsequence,\vspace*{1pt} $\int_0^tf(X_{s-}) \, d\overline
{A}{}^{d,\ell_k}_s$
[resp., $\int_0^tf(X_{s-}) \, d\overline{A}{}^{d,\ell_k}_s$]
uniformly converges to $\int_0^tf(X_s)\circ d\overline{A}{}^{d}_s$
[resp., $\int_0^tf(X_s)\circ d\overline{A}{}^{d}_s$]
on each compact subinterval of $[0,\infty[$ ${\mathbb{P}}_x$-a.s. for
q.e. $x\in E$.
So the conclusion holds by Corollary \ref{cor:compensator3}.
For general $M\in(\mathcal{M}_{\mathrm{loc}})_*$ [resp., $M\in(\mathcal
{M}^{\auf 0,\zeta\auf }_{\mathrm{loc}})_*$]
with its jump function $\varphi\in\mathcal{J}_*$ (resp., $\varphi\in
\widehat{\mathcal{J}}_*$), the assertion
is clear from Lemma \ref{lem:localNakaoInt}.
\end{pf}


Recalling
Theorem \ref{thm:purediscontimartincorrespondince} and the last
description of
Section \ref{sec:preliminary},
for $u\in\dot{\mathcal{F}}^{\ddag}_{\mathrm{loc}}$
(resp., $u\in\dot{\mathcal{F}}_{\mathrm{loc}}^{\dag}$),
there exists an $M^{u,d}\in\mathcal{M}^{d}_{\mathrm{loc}}(\subset
\accentset{\circ}{\mathcal{M}}^d_{f\mbox{-}\mathrm{loc}})$
[resp., $M^{u,d}\in\mathcal{M}^{d,\auf 0,\zeta\auf }_{\mathrm{loc}}(
\subset\accentset{\circ}{\mathcal{M}}^d_{\mathrm{loc}})$]
such that
$M^{u,d}_t-M^{u,d}_{t-}=u(X_t)-u(X_{t-})$ for all $t\in[0,\infty[$
(resp., $[0,\zeta[$) ${\mathbb{P}}_x$-a.s. for
q.e. $x\in E$.
By Lemma \ref{lem:MAFcontpartvanish}(1),
we can define $M^{u,c}\in\mathcal{M}^c_{\mathrm{loc}}=\accentset{\circ}{\mathcal{M}}^c_{f\mbox{-}\mathrm{loc}}$
for $u\in\dot{\mathcal{F}}_{\mathrm{loc}}$;
$M_t^{u,c}:=M_t^{u_n,c}$ for
$t<\sigma_{E\setminus G_n}$ for some $\{G_n\}\in\Theta$
and $u_n\in\mathcal{F}$ such that $u=u_n$ $m$-a.e.  on $G_n$
for each $n\in\mathbb{N}$ (see Remark \ref{rem:coincidenceMAFs}).
Put $A^{u,c}:=M^{u,c}+\Gamma(M^{u,c})$ for $M^{u,c}$, which can be
defined by
way of Lemma \ref{lem:divlikeCAFlocal}, and $A^{u,d}:=M^{u,d}+\Gamma
(M^{u,d})$ for $M^{u,d}$, which is defined by
Definition \ref{df:localizationGamma}.
\begin{cor}\label{cor:purejumpapprox2}
Take $f\in\dot{\mathcal{F}}_{\mathrm{loc}}$ and $u\in\dot{\mathcal
{F}}^{\ddag}_{\mathrm{loc}}$.
Suppose that $f$ is extended to be a real-valued function on
$E_{\partial}$
with $f(\partial):=0$.
Then
%
\begin{eqnarray}
\label{eq:GeneFukStra}
\int_0^tf(X_{s-})  \,dA_s^{u,d}&=&
{\sum_{0<s\leq t}}^{\!\!\!\!*}f(X_{s-})
\bigl(u(X_s)-u(X_{s-})\bigr), \\
\label{eq:GeneFukStra*}
\int_0^tf(X_s)\circ dA_s^{u,d}&=&
{\sum_{0<s\leq t}}^{\!\!\!\!*}\frac{f(X_s)+f(X_{s-})}{2}
\bigl(u(X_s)-u(X_{s-})\bigr)
\end{eqnarray}
hold for all $t\in[0,\infty[$ ${\mathbb{P}}_x$-a.s. for q.e. $x\in E$.
Similarly
%
\begin{equation} \label{eq:GEneFukuDec*}
A_t^{u,d}=
{\sum_{0<s\leq t}}^{\!\!\!\!*}\bigl(u(X_s)-u(X_{s-})\bigr)
\end{equation}
hold for all $t\in[0,\infty[$ ${\mathbb{P}}_x$-a.s. for q.e. $x\in E$.
More generally, if $u\in\dot{\mathcal{F}}_{\mathrm{loc}}^{\dag}$ and $f\in
\dot{\mathcal{F}
}_{\mathrm{loc}}$ is
only defined on $E$, then all assertions above hold for all
$t\in[0,\zeta[$ ${\mathbb{P}}_x$-a.s. for q.e. $x\in E$.
\end{cor}

Owing to (\ref{eq:GEneFukuDec*}), we can obtain a generalized
Fukushima decomposition for $u\in\dot{\mathcal{F}}^{\dag}_{\mathrm{loc}}$:
\begin{theorem}[(Generalized Fukushima decomposition)]\label{thm:FDlocallyinDir}
For $u\in\dot{\mathcal{F}}^{\ddag}_{\mathrm{loc}}$,
the additive functional $A^u$ defined by $A_t^u:=u(X_t)-u(X_0)$ can be
decomposed as
\[
A^u=M^u+N^u, \qquad M^u\in\mathcal{M}_{\mathrm{loc}},\qquad  N^u\in\mathcal
{N}_{c,f\mbox{-}\mathrm{loc}}
\]
in the sense that $A_t^u=M_t^u+N_t^u,  t\in[0,\infty[\mbox{ }  {\mathbb
{P}}_x\mbox
{-a.s. for q.e. }x\in E$.
More generally, if $u\in\dot{\mathcal{F}}_{\mathrm{loc}}^{\dag}$, then $A^u$ is
decomposed as
\[
A^u=M^u+N^u,\qquad  M^u\in\mathcal{M}^{\auf 0,\zeta\auf }_{\mathrm{loc}},\qquad  N^u\in
\mathcal
{N}_{c,\mathrm{loc}}
\]
in the sense that $A_t^u=M_t^u+N_t^u,  t\in[0,\zeta[ \mbox{ } {\mathbb
{P}}_x\mbox
{-a.s. for q.e. }x\in E$.
Such decompositions are unique up to the equivalence of (local)
additive functionals.
\end{theorem}
\begin{pf}
The uniqueness is proved in Corollary \ref{cor:uniqeness}.
We shall only prove the existence in the first assertion.
We set $M^{u}:=M^{u,c}+M^{u,d}\in(\mathcal{M}_{\mathrm{loc}})_{\mathrm{as}}$ and
$N^u:=\Gamma(M^u)\in\mathcal{N}_{c,f\mbox{-}\mathrm{loc}}$, where
$M^{u,c}$ and $M^{u,d}$ are defined above.
Take $\{G_n\}\in\Theta$ and
$\{u_n\}\subset\mathcal{F}_b$ such that $u-u(\partial)=u_n$
$m$-a.e. on $G_n$.
Then (\ref{eq:GEneFukuDec*}) implies that for $t\in[0,\sigma
_{E\setminus G_n}[$
\begin{eqnarray*}
u(X_t)-u(X_0)&=&u_n(X_t)-u_n(X_0)=A_t^{u_n,c}+A_t^{u_n,d}\\
&=&
A_t^{u_n,c}+{\sum_{0<s\leq t}}^{\!\!\!\!*}\bigl(u_n(X_s)-u_n(X_{s-})\bigr)\\
&=&A_t^{u,c}+{\sum_{0<s\leq t}}^{\!\!\!\!*}\bigl(u(X_s)-u(X_{s-})\bigr)
\stackrel{\mbox{\fontsize{8.36}{10.36}\selectfont{(\ref{eq:GEneFukuDec*})}}}{=} A_t^{u,c}+A_t^{u,d}\\
&=&M_t^{u}+\Gamma(M^{u})_t
\end{eqnarray*}
${\mathbb{P}}_x$-a.s. for q.e. $x\in E$.
\end{pf}
\begin{rem}\label{rem:GEneFuku}
(1)
We emphasize that $\mathbf{1}_{E_{\partial}}$ does not satisfy
$\mathbf{1}_{E_{\partial}}(\partial)=0$.
So we

{\smallskipamount=0pt
\begin{enumerate}[(2)]
\item[ ] cannot deduce (\ref{eq:GEneFukuDec*}) from
(\ref{eq:GeneFukStra}), (\ref{eq:GeneFukStra*}).
\item[(2)] For $f\in\dot{\mathcal{F}}_{\mathrm{loc}}$ with $f(\partial)=0$
and $u\notin\dot{\mathcal{F}}_{\mathrm{loc}}^{\dag}$, we have no way to define
$M^{u,d}$, $\Gamma(M^{u,d})$ and stochastic integrals with respect to them.
However, we can define the left-hand sides of
(\ref{eq:GeneFukStra}) and (\ref{eq:GeneFukStra*}) keeping the same
expressions as they have.
\item[(3)] In Theorem \ref{thm:FDlocallyinDir},
$M^u$ for $u\in\dot{\mathcal{F}}_{\mathrm{loc}}^{\ddag}$ (resp., $u\in\dot
{\mathcal{F}
}_{\mathrm{loc}}^{\dag}$) can be
decomposed to $M_t^u=M_t^{u,c}+M_t^{u,j}+M_t^{u,\kappa}$ $t\in
[0,\infty[$
(resp., $M_t^u=M_t^{u,c}+M_t^{u,d}$ $t\in[0,\zeta[$)
${\mathbb{P}}_x$-a.s. for q.e. $x\in E$, where $M^{u,i}\in\mathcal
{M}_{\mathrm{loc}}^i$,
$i=c,j,\kappa$
(resp., $M^{u,c}\in\mathcal{M}_{\mathrm{loc}}^{c,\auf 0,\zeta\auf }$,
$M^{u,d}\in\mathcal{M}
_{\mathrm{loc}}^{d,\auf 0,\zeta\auf }$).
\end{enumerate}}
\end{rem}

Now we expose a generalized It\^o formula in terms of our stochastic integrals.
\begin{theorem}[(Generalized It\^o formula)]\label{thm:ItofornulaFS}
Suppose that $\Phi\in C^1 (\mathbb{R}^N)$ and take
$u=( u_1, \ldots, u_N) \in(\dot{\mathcal{F}}_{\mathrm{loc}})^N$ having an
$\mathbb{R}
^N$-valued extension on $E_{\partial}$.
Then:
\begin{enumerate}[(1)]
\item[(1)]
$\Phi(u)\in\dot{\mathcal
{F}}_{\mathrm{loc}}$ and
for each $k=1,2,\ldots,N$,
$\Phi_k(u)\in L^2_{\mathrm{loc}}(\{G_n\};\mu_{\langle M^{u_k,c}\rangle})$ for
some $\{G_n\}\in\Theta$, where
$\Phi_k:=\frac{\partial\Phi}{\partial x_k}$ and
%
%
\begin{equation}\label{eqn:ItocontiStratonovich}
A_t^{\Phi(u),c}=\sum_{k=1}^N\int_0^t\Phi_k(u(X_s))\circ dA_s^{u_k,c}
\end{equation}
holds for all $t\in[0,\infty[$ ${\mathbb{P}}_x$-a.s. for q.e. $x\in
E$. If we
further assume $\Phi\in C^2(\mathbb{R}^N)$, then for each $k,\ell
=1,2,\ldots,N$,
$\Phi_k(u)\in\dot{\mathcal{F}}_{\mathrm{loc}}$, $\Phi_{k\ell}(u)\in
L^2_{\mathrm{loc}}(\{G_n\};\mu_{\langle M^{u_k,c}\rangle})$ for some $\{G_n\}
\in\Theta$, where
$\Phi_{k\ell}:=\frac{\partial^2\Phi}{\partial x_k\partial x_{\ell
}}$, and
%
\begin{eqnarray}\label{eqn:Itoconti}
A_t^{\Phi(u),c}&=&
\sum_{k=1}^N\int_0^t\Phi_k(u(X_{s-})) \,
dA_s^{u_k,c}\nonumber\\[-8pt]\\[-8pt]
&&{} +
\frac{1}{2}\sum_{k,\ell=1}^N\int_0^t\Phi_{k\ell}(u(X_s))\,d\langle
M^{u_k,c},M^{u_{\ell},c}\rangle_s
\nonumber
\end{eqnarray}
holds for all $t\in[0,\infty[$ ${\mathbb{P}}_x$-a.s. for q.e. $x\in E$.
\item[(2)]
Suppose $u\in(\dot{\mathcal{F}
}_{\mathrm{loc}}^{\ddag})^N$ and
$\Phi\in C^2(\mathbb{R}^N)$. Then $\Phi_k(u)\in\dot{\mathcal
{F}}_{\mathrm{loc}}$ for each
$k=1,2,\ldots,N$.
Moreover, if we assume
$\Phi(u)\in\dot{\mathcal{F}}_{\mathrm{loc}}^{\ddag}$, then
%
\begin{eqnarray}
\label{eqn:ItojumpIto}
A_t^{\Phi(u),d}&=&\sum_{k=1}^N\int_0^t\Phi_k(u(X_{s-}))
\,dA_s^{u_k,d}\nonumber\\[-8pt]\\[-8pt]
&&{}
+ \sum_{s\leq t} \Biggl( \Delta\Phi(u(X_s))
-\sum_{k=1}^N \Phi_k
(u(X_{s-})) \Delta u_k(X_s) \Biggr)\nonumber
\end{eqnarray}
and
%
\begin{eqnarray}\label{eqn:ItojumpStratonoovich}
A_t^{\Phi(u),d}&=&\sum_{k=1}^N\int_0^t\Phi_k(u(X_s))\circ
\,dA_s^{u_k,d}\nonumber\\[-8pt]\\[-8pt]
&&{}
+ \sum_{s\leq t} \Biggl( \Delta\Phi(u(X_s))
-\sum_{k=1}^N \frac{\Phi_k(u(X_s))+ \Phi_k
(u(X_{s-}))}{2} \Delta u_k(X_s) \Biggr)\hspace*{-15pt}\nonumber
\end{eqnarray}
hold for all $t\in[0,\infty[$ ${\mathbb{P}}_x$-a.s. for q.e. $x\in E$,
where the last terms in the right-hand sides are absolutely
convergent uniformly on each compact interval of $[0,\infty[$
${\mathbb{P}}_x$-a.s. for q.e. $x\in E$. If we replace $\dot{\mathcal
{F}}_{\mathrm{loc}}^{\ddag
}$ with $\dot{\mathcal{F}}_{\mathrm{loc}}^{\dag}$
in the above conditions,
then formulas (\ref{eqn:ItojumpIto}) and (\ref
{eqn:ItojumpStratonoovich}) hold only on
$[0,\zeta[$ ${\mathbb{P}}_x$-a.s. for q.e. $x\in E$ without assuming
the $\mathbb{R}
^N$-valued extension of $u$
on $E_{\partial}$.
\item[(3)]
Under the same conditions in (2), we have
$\Phi_k(u)\in\dot{\mathcal{F}}_{\mathrm{loc}}$ for each $k=1,2,\ldots,N$,
%
\begin{eqnarray}\label{eqn:ItoIto}
A_t^{\Phi(u)}&=&\sum_{k=1}^N \int_0^t \Phi_k( u(X_{s-}) )\,
dA^{u_k}_s
\nonumber\\
&&{}
+\frac{1}{2}\sum_{k,\ell=1}^N\int_0^t\Phi_{k\ell}(u(X_s))\,d\langle
M^{u_k,c},M^{u_{\ell},c}\rangle_s
\\
&&{}+ \sum_{s\leq t} \Biggl( \Delta\Phi(u(X_s))
-\sum_{k=1}^N \Phi_k
(u(X_{s-}))\Delta u_k(X_s) \Biggr)\nonumber
\end{eqnarray}
and
%
\begin{eqnarray} \label{eqn:ItoStratonovich}
A_t^{\Phi(u)}&=&
\sum_{k=1}^N \int_0^t \Phi_k( u(X_{s}) ) \circ
dA^{u_k}_s \nonumber\\[-8pt]\\[-8pt]
&&{}
+ \sum_{s\leq t} \Biggl( \Delta\Phi(u(X_s))
-\sum_{k=1}^N \frac{\Phi_k(u(X_s))+ \Phi_k
(u(X_{s-}))}{2} \Delta u_k(X_s) \Biggr) \hspace*{-15pt}\nonumber
\end{eqnarray}
hold for all $t\in[0,\infty[$ ${\mathbb{P}}_x$-a.s. for q.e. $x\in E$,
where the last terms in the right-hand sides
are absolutely convergent uniformly on each compact interval of
$[0,\infty[$
${\mathbb{P}}_x$-a.s. for q.e. $x\in E$. If we replace $\dot{\mathcal
{F}}_{\mathrm{loc}}^{\ddag
}$ with $\dot{\mathcal{F}}_{\mathrm{loc}}^{\dag}$
in the above conditions,
then formulas (\ref{eqn:ItoIto}) and (\ref{eqn:ItoStratonovich}) hold
only on
$[0,\zeta[$ ${\mathbb{P}}_x$-a.s. for q.e. $x\in E$ without assuming
the $\mathbb{R}^N$-valued
extension of $u$ on $E_{\partial}$.
\end{enumerate}
\end{theorem}

We call (\ref{eqn:ItoIto}) the \textit{It\^o formula for It\^o-type
integrals} and
(\ref{eqn:ItoStratonovich}) the \textit{It\^o formula for
Fisk--Stratonovich-type integrals}.
\begin{cor}[(Chain and Leibniz rules for purely discontinuous
part)]\label{cor:chainrulemartingale}
Under the same conditions as in Theorem \ref{thm:ItofornulaFS}(2), we have
that
%
\begin{eqnarray}\qquad
\label{eq:chainrulemartingale}
M_t^{\Phi(u),d}&=&\sum_{k=1}^N \int_0^t \Phi_k( u(X_{s-}) )
\,dM^{u_k,d}_s \nonumber\\
&&{}+\sum_{s\leq t} \Biggl( \Delta\Phi(u(X_s))
-\sum_{k=1}^N \Phi_k
(u(X_{s-}))\Delta u_k(X_s) \Biggr)\\
&&{}- \Biggl\{\sum_{s\leq t} \Biggl( \Delta\Phi(u(X_s))
-\sum_{k=1}^N \Phi_k
(u(X_{s-}))\Delta u_k(X_s) \Biggr) \Biggr\}^p,\nonumber
\\
\label{eq:chainrulemartingale*}
M_t^{\Phi(u),d}&=&\sum_{k=1}^N \int_0^t \Phi_k( u(X_{s}) )\circ
\,dM^{u_k,d}_s \nonumber\\
&&{}
+\sum_{s\leq t} \Biggl( \Delta\Phi(u(X_s))\nonumber\\[-8pt]\\[-8pt]
&&{}
\hspace*{31.7pt}-\sum_{k=1}^N \frac{\Phi_k(u(X_s))+\Phi_k
(u(X_{s-}))}{2}\Delta u_k(X_s) \Biggr)\nonumber\\
&&{}
- \Biggl\{\sum_{s\leq t} \Biggl( \Delta\Phi(u(X_s))
-\sum_{k=1}^N \Phi_k
(u(X_{s-}))\Delta u_k(X_s) \Biggr) \Biggr\}^p\nonumber
\end{eqnarray}
and
%
\begin{eqnarray}\label{eq:chainrulemartingale**}\qquad
\Gamma\bigl(M^{\Phi(u),d}\bigr)_t&=&\sum_{k=1}^N \int_0^t \Phi_k( u(X_{s}) )
\,d\Gamma(M^{u_k,d})_s
\nonumber\\[-8pt]\\[-8pt]
&&{}+ \Biggl\{\sum_{s\leq t} \Biggl( \Delta\Phi(u(X_s))
-\sum_{k=1}^N \Phi_k
(u(X_{s-}))\Delta u_k(X_s) \Biggr) \Biggr\}^p \nonumber
\end{eqnarray}
hold for all $t\in[0,\infty[$ ${\mathbb{P}}_x$-a.s. for q.e. $x\in E$,
where $B_t^p$ denotes the dual predictable projection of $B_t$ for an
AF $B$.
If we replace $\dot{\mathcal{F}}_{\mathrm{loc}}^{\ddag}$ with $\dot{\mathcal
{F}}_{\mathrm{loc}}^{\dag}$
in the conditions, (\ref{eq:chainrulemartingale}) and (\ref
{eq:chainrulemartingale*}) hold
only on
$[0,\zeta[$ ${\mathbb{P}}_x$-a.s. for q.e. $x\in E$ without assuming
the $\mathbb{R}^N$-valued
extension of $u$ on $E_{\partial}$.
In particular, for $u,v\in\dot{\mathcal{F}}_{\mathrm{loc}}^{\dag}$
(resp., $u,v\in
\dot{\mathcal{F}}_{\mathrm{loc}}^{\ddag}$) with $uv\in
\dot{\mathcal{F}}_{\mathrm{loc}}^{\dag}$ [resp., $uv\in\dot{\mathcal
{F}}_{\mathrm{loc}}^{\ddag}$ with
$u(\partial)=v(\partial)=0$],
%
\begin{eqnarray}\label{eq:Leibnizrulemartingale}
M_t^{uv,d}&=&\int_0^tu(X_{s-})  \,dM^{v,d}_s+\int_0^tv(X_{s-})  \,dM^{u,d}_s
\\
&&{}+[M^{u,d},M^{v,d}]_t-\langle M^{u,d},M^{v,d}\rangle_t
\nonumber\\
\label{eq:Leibnizrulemartingale*}
&=&\int_0^tu(X_{s-})\circ dM^{v,d}_s+\int
_0^tv(X_{s-})\circ dM^{u,d}_s\nonumber\\[-8pt]\\[-8pt]
&&{}-
\langle M^{u,d},M^{v,d}\rangle_t,\nonumber
\\
\label{eq:Leibnizrulemartingale**}
\Gamma(M^{uv,d})_t&=&\int_0^tu(X_{s})  \,d\Gamma(M^{v,d})_s+
\int_0^tv(X_{s})  \,d\Gamma(M^{u,d})_s\\
&&{}
+\langle M^{u,d},M^{v,d}\rangle_t
\nonumber
\end{eqnarray}
hold for all $t\in[0,\zeta[$ (resp., $t\in[0,\infty[$) ${\mathbb{P}}
_x$-a.s. for q.e. $x\in E$.
\end{cor}
\begin{cor}[(Fisk--Stratonovich integration by parts formula)]
\label{cor:integtrationbyparts}
For $u,v\in\dot{\mathcal{F}}_{\mathrm{loc}}^{\dag}$ (resp., $u,v\in\dot
{\mathcal{F}
}_{\mathrm{loc}}^{\ddag}$) with $uv\in
\dot{\mathcal{F}}_{\mathrm{loc}}^{\dag}$ [resp., $uv\in\dot{\mathcal
{F}}_{\mathrm{loc}}^{\ddag}$ with
$u(\partial)=v(\partial)=0$],
%
\begin{equation}\label{eq:integtrationbyparts}\quad
u(X_t)v(X_t)-u(X_0)v(X_0)=\int_0^tu(X_s)\circ dA^v_s+\int
_0^tv(X_s)\circ dA^u_s
\end{equation}
holds for all $t\in[0,  \zeta[$ (resp., $t\in[0,\infty[$) ${\mathbb{P}}
_x$-a.s. for q.e. $x\in E$.
\end{cor}
\begin{rem}\label{rem:ItofornulaFS}
(1)
In \cite{CFKZStoch},
we prove a generalized It\^o formula for
$u\in(\mathcal{F})^N$

{\smallskipamount=0pt
\begin{enumerate}[(2)]
\item[ ] under the law for
quasi everywhere starting points, extending the early result by
Nakao \cite{Na}. Our It\^o formula can be applied to a wider class of
integrators than that in \cite{CFKZStoch}.
\item[(2)]
In Theorem \ref{thm:ItofornulaFS}(2),
if $u\in((\dot{\mathcal{F}}_{\mathrm{loc}})_b)^N$ with $u(\partial)\in
\mathbb{R}^N$,
$\Phi\in C^2(\mathbb{R}^N)$, then $\Phi(u)\in\dot{\mathcal
{F}}_{\mathrm{loc}}^{\ddag}$.\vspace*{1pt}
\item[(3)]
$\Phi_k(u)\in\dot{\mathcal{F}}_{\mathrm{loc}}$ does not necessarily satisfy
$\Phi_k(u)(\partial)=0$. So the integrals appeared in the
first terms of the right-hand sides of (\ref{eqn:Itoconti})--(\ref
{eqn:ItoIto}) should be understood
to be modified, for example,
$\int_0^t\Phi_k(u(X_s))\,dA_s^{u_k}$ should be understood as
$\int_0^t(\Phi_k(u(X_s))-\Phi_k(u(\partial)))\,dA_s^{u_k}+\Phi
_k(u(\partial))A_t^{u_k}$.
\item[(4)]
Comparing with (\ref{eqn:ItocontiStratonovich}), the case for
diffusion part,
our It\^o formulas, (\ref{eqn:ItojumpStratonoovich}) and (\ref
{eqn:ItoStratonovich}), for
Fisk--Stratonovich integrals are not so simple.
This phenomenon
can be found in the It\^o formula for Fisk--Stratonovich integral
exposed in
II. Theorem 34 and V. Theorem 21 of \cite{ProtStoch} in the framework
of semi-martingales.
We emphasize that the expression of the
second term (denoted by $C_t$) of the right-hand side in (\ref
{eqn:ItoStratonovich}) is different
from theirs [i.e., the
third term of the right-hand side in (\ref{eqn:ItoIto}), which is
the usual expression of the It\^o formula for purely discontinuous part].
Note that $C_t$
is an \textit{odd additive functional},
that is, for each $t>0$, $C_t\circ r_t+C_t=0$ ${\mathbb{P}}_m$-a.e. on
$\{t<\zeta\}$, where
$r_t$ is the time reverse operator.
Hence, both sides in our formula (\ref{eqn:ItoStratonovich})
possess this property, which is not yielded by the It\^o
formula in \cite{ProtStoch}.

\item[(5)]
In Theorems \ref{thm:Stratonovichlocalrepre}, \ref{thm:ItofornulaFS} and
Corollaries \ref{cor:purejumpapprox2}--\ref{cor:integtrationbyparts},
we do not require
the strict $\mathcal{E}$-quasi-continuities of $f$, $u$ and $v$.
If we do not impose the condition that such functions are extended on
$E_{\partial}$
and vanish on $\{\partial\}$, or if we only assume that $(\mathcal
{E},\mathcal{F})$ is
not necessarily regular
(i.e., quasi-regularity only holds), then all assertions are restricted
to ``for all
$t\in[0,\zeta[$ ${\mathbb{P}}_x$-a.s. for q.e. $x\in E$'' and
each convergence
of the right-hand side is uniform
on compact subinterval of $[0,\zeta[$ ${\mathbb{P}}_x$-a.s. for
q.e. $x\in E$.

\item[(6)]
In \cite{FolIto}, an It\^o formula for a general multi-dimensional
process with finite
quadratic variation is presented, but the formula like
(\ref{eqn:ItoStratonovich}) is not exposed in~\cite{FolIto}.

\item[(7)]
As noted in Remark \ref{rem:GEneFuku},
even for $u\notin(\dot{\mathcal{F}}_{\mathrm{loc}}^{\dag})^N$, we can define the
first terms in the right-hand sides
of (\ref{eqn:ItojumpIto}), (\ref{eqn:ItojumpStratonoovich}), (\ref
{eqn:ItoIto}) and
(\ref{eqn:ItoStratonovich}). So the formulas hold in this setting
without using
stochastic integrals with respect to $M^{u,d}$, $\Gamma(M^{u,d})$.
So the conclusion of Corollary \ref{cor:integtrationbyparts} also
holds for $u,v\in\dot{\mathcal{F}}_{\mathrm{loc}}$ with
$u(\partial)=v(\partial)=0$ in this context.
\end{enumerate}}
\end{rem}
\begin{pf*}{Proof of Theorem \ref{thm:ItofornulaFS}}
(2) is a consequence of (1) and (3). We first prove~(1).
The former assertion of (1) follows from
Theorems 6.1 and 7.2 in
\cite{Kwfunc}. Note that $M^{\Phi(u),c}, M^{u_k,c}\in\mathcal{M}
_{\mathrm{loc}}^c=\accentset{\circ}{\mathcal{M}}^c_{f\mbox{-}\mathrm{loc}}$ and
$\Gamma(M^{\Phi(u),c}),\Gamma(M^{u_k,c})\in\mathcal{N}^c_{f\mbox
{-}\mathrm{loc}}$, which are defined on $[0,\infty[$ under
${\mathbb{P}}_x$ for q.e. $x\in E$.
Formula (\ref{eqn:ItocontiStratonovich}) can be obtained from the
chain rule for
continuous part of MAF
%
\begin{equation}\label{eq:chainrulecontmartingale}
M_t^{\Phi(u),c}=\sum_{k=1}^N\int_0^t\Phi_k(u(X_{s-}))\,dM_s^{u_k,c}
\end{equation}
for all $t\in[0,\infty[
{\mathbb{P}}_x\mbox{-a.s. for q.e. }x\in E$
(see Theorem 7.2 in \cite{Kwfunc}).
The latter assertion of (1) also follows from
Theorem 7.2 in \cite{Kwfunc} and (\ref{eq:contnakaoint}).

Next we prove (3). Applying Theorem 6.1 in
\cite{Kwfunc} to $\Phi_k\in C^1(\mathbb{R}^N)$ again,
we have $\Phi_k(u)\in\dot{\mathcal{F}}_{\mathrm{loc}}$ for $u\in(\dot
{\mathcal{F}
}_{\mathrm{loc}}^{\dag})^N$. Equation
(\ref{eqn:ItoIto}) is proved by Nakao \cite{Na} for the case $u\in
(\mathcal{F}_b)^N$.
(\ref{eqn:ItoStratonovich}) for $u\in(\mathcal{F}_b)^N$ also follows from
(\ref{eqn:ItoIto}) for
$u\in(\mathcal{F}_b)^N$ and
that for each $k=1,2,\ldots,N$
\[
\int_0^t\Phi_k(u(X_s))\circ dM_s^{u_k,d}=\int_0^t\Phi
_k(u(X_{s-}))\,dM_s^{u_k,d}+\frac{1}{2}
\bigl[M^{\Phi_k(u),d}, M^{u_k,d}\bigr]_t
\]
for all $t\in[0,\infty[$ ${\mathbb{P}}_x$-a.s. for q.e. $x\in E$.
Equations (\ref{eqn:ItoIto}) and
(\ref{eqn:ItoStratonovich}) for general $u\in\dot{\mathcal
{F}}_{\mathrm{loc}}^{\ddag}$
(or $u\in\dot{\mathcal{F}}_{\mathrm{loc}}^{\dag}$) hold for all $t\in
[0,\infty[$
(resp., $t\in[0,\zeta[$) ${\mathbb{P}}_x$-a.s. for
q.e. $x\in E$
in view of the both of local properties of stochastic integrals by
$M^{u_k,d}$ and $\Gamma(M^{u_k,d})$ (see Lemma \ref{lem:localNakaoInt}).
\end{pf*}
\begin{ex}[(Symmetric L\'evy process on $\mathbb{R}^N$)]
\hspace*{-1.7pt}Let $\mathbf{M}=(\Omega,X_t,\zeta, {\mathbb{P}}_x)\hspace*{-0.5pt}_{x\in
\mathbb{R}^N}$ be the
symmetric L\'evy process. That is,
\textbf{M} is a time homogeneous additive process determined by a family
$\{\nu_t\}$ of probability measures on $\mathbb{R}^N$ satisfying (4.17),
(4.18) and (4.19) in \cite{FOT}.
Let $(\mathcal{E},\mathcal{F})$ be the corresponding Dirichlet form
on $L^2(\mathbb{R}^N)$.
Then
$(\mathcal{E},\mathcal{F})$ is given by
\begin{eqnarray*}
\cases{
\displaystyle\mathcal{F}= \biggl\{u\in L^2(\mathbb{R}^N) \Bigm|  \int_{\mathbb{R}^N}
| \hat{u}(\xi)|^2\psi(\xi)\,d\xi<\infty
\biggr\},\cr
\displaystyle\mathcal{E}(u,v)=\int_{\mathbb{R}^N}\hat{u}(\xi)\overline{\hat
{v}}(\xi)\psi
(\xi)\,d\xi, \qquad  u,v\in\mathcal{F},}
\end{eqnarray*}
where $\hat{u}(\xi):=\frac1{(2\pi)^{d/2}}\int_{\mathbb
{R}^N}e^{i\langle
\xi,x\rangle}u(x)\,dx$ and
$\psi(x)$ is the function determined by ${\mathbb{E}}_0[e^{i\langle
\xi
,X_t\rangle}]=e^{-t\psi(\xi)}$.
We assume that \textbf{M} is purely discontinuous; namely $\psi$ has the
following expression:
\[
\psi(\xi)=\int_{\mathbb{R}^N}(1-\cos\langle\xi,\eta\rangle)\nu
(d\eta),
\]
where $\nu$ is a symmetric measure on $\mathbb{R}^N\setminus\{0\}$
such that
$\int_{\mathbb{R}^N\setminus\{0\}} (|\xi|^2\wedge1)\nu(d\xi
)<\infty$,
which is called
the \textit{L\'evy measure} of \textbf{M}.
We see $C_0^{\mathrm{Lip}}(\mathbb{R}^N)\subset\mathcal
{F}$, hence
$C^{\mathrm{Lip}}_{\mathrm{loc}}(\mathbb{R}^N)\subset\mathcal
{F}_{\mathrm{loc}}$, because,
in view of Corollary 7.16 in \cite{BFPotential},
\[
1+\psi(\xi)\leq c(1+|\xi|^2)  \qquad{\forall}\xi\in\mathbb{R}^N
\]
for some constant $c>0$. Here $C_0^{\mathrm{Lip}}(\mathbb{R}^N)$
[resp., $C_{\mathrm{loc}}^{\mathrm{Lip}}(\mathbb{R}^N)$] is the family of
Lipschitz continuous functions with compact support
(resp., locally Lipschitz continuous functions) and $\mathcal{F}_{\mathrm{loc}}$
is the
space of functions locally in $\mathcal{F}$
in the ordinary sense (see \cite{FOT}).
Further
$(\mathcal{E},\mathcal{F})$ is a regular Dirichlet form having
$C_0^{\infty}(\mathbb{R}^N)$
as its core
(see \cite{Takexit}).
Define $N(x,A):=\nu(A-x)$, $N(x,\{\partial\})=0$ for $A\in\mathscr
{B}(\mathbb{R}
^N)$, $x\in\mathbb{R}^N$
and $H_t=t$.
By Theorem 19.2(i) in \cite{SatoLevy}, we have
\[
N(x,A)={\mathbb{E}}_x \biggl[\sum_{0<s\leq1}\mathbf{1}_{A}(X_s-X_{s-})
\biggr],\qquad
A\in\mathscr{B}
(\mathbb{R}^N),
\]
and hence $(N,H)$ becomes a L\'evy system of \textbf{M} (see also Section
7 in \cite{KWsquare}).
By Corollary \ref{cor:purejumpapprox2}, we have that for any
$u\in C_{\mathrm{loc}}^{\mathrm{Lip}}(\mathbb{R}^N)$
\[
u(X_t)-u(X_0)=
{\sum_{0<s\leq t}}^{\!\!\!\!*}\bigl(u(X_s)-u(X_{s-})\bigr)
\]
holds for all $t\in[0,\zeta[$ ${\mathbb{P}}_x$-a.s. for q.e. $x\in
\mathbb{R}^N$.
Further we assume $\nu(dy)=f(|y|)\,dy$, where $f$ is a
Borel function satisfying (\ref{eq:conditionjump}).
Let $u\in C^{\mathrm{Lip}}(\mathbb{R}^N)$. Then
\[
\sup_{x\in K}\int_{\mathbb{R}^N}\bigl(u(x+y)-u(x)\bigr)^2\nu(dy)<\infty
\qquad\mbox{for
any compact set $K$,}
\]
hence $u\in\mathcal{F}_{\mathrm{loc}}^{\dag}$. Therefore $u$ admits the generalized
Fukushima decomposition.
\end{ex}
\begin{ex}[(Symmetric stable process on $\mathbb{R}^N$)]
We fix $\alpha\in\ ]0,2[$.
Let $\mathbf{M}=(\Omega, X_t, {\mathbb{P}}_x)_{x\in\mathbb{R}^N}$
be a
L\'evy process on $\mathbb{R}^N$ with
\[
{\mathbb{E}}_{0}\bigl[e^{\sqrt{-1}\langle\xi, X_{t} \rangle}\bigr]=
e^{-t|\xi|^{\alpha}}.
\]
\textbf{M} is called the
\textit{symmetric $\alpha$-stable process}.
It is known that \textbf{M} is conservative.
Let $(\mathcal{E},\mathcal{F})$ be the associated Dirichlet form on
$L^2(\mathbb{R}^N)$ with
\textbf{M},
which is given by
%
\begin{equation}\label{eq:stabledirichlet}\qquad\quad
\cases{\displaystyle\mathcal{F}= \biggl\{u\in L^2(\mathbb{R}^N)  \Bigm|
\iint_{\mathbb{R}^N\times
\mathbb{R}^N}
\frac{(u(x)-u(y))^2}{\vert x-y\vert^{N+\alpha}}\,dx\,dy<\infty
\biggr\},\vspace*{2pt}\cr
\displaystyle\mathcal{E}(u,v)= \frac{A(N,-\alpha)}{2}\iint_{\mathbb
{R}^N\times\mathbb{R}^N}
\frac{(u(x)-u(y))(v(x)-v(y))}{\vert x-y\vert^{N+\alpha}}\,dx\,dy, \cr
\qquad u,v\in\mathcal{F},}
\end{equation}
where
$
A(N,\gamma):=\frac{|\gamma|\Gamma(({N-\gamma})/{2})}
{2^{1+\gamma}\pi^{N/2}\Gamma(1+{\gamma}/{2})}, \gamma<N.
$
The L\'evy system $(N,H)$ of \textbf{M} is given by $N(x,dy):=A (N,-\alpha
)\vert
x-y\vert^{-(N+\alpha)}\,dy$ and $H_t=t$. So $\mu_H(dx)=dx$.
Hence $J(dx\,dy)=f(|x-y|)\,dx\,dy$ for $f(r):=A (N,-\alpha) r^{-N-\alpha}$, $r>0$.
Note that $f$ does not satisfy (\ref{eq:conditionjump}).
Take $\beta\in[0,\alpha[$. Assume that $N\geq\alpha$, hence $\{0\}
$ is polar, and take $u\in C^1(\mathbb{R}^N\setminus\{0\})\cap
C^{0,\beta
/2}(\mathbb{R}^N)$. Here $C^{0,\beta/2}(\mathbb{R}^N)$ is the family
of $\beta
/2$-H\"older continuous functions on $\mathbb{R}^N$.
For example, for a function $F\in C^1([0,\infty[)$ with bounded
derivative $F'$, $u(x):=F(|x|^{\beta/2})$ is a function in
$C^1(\mathbb{R}^N\setminus\{0\})\cap C^{0,\beta/2}(\mathbb{R}^N)$.
Then $u\in\dot{
\mathcal{F}}_{\mathrm{loc}}^{\dag}=\dot{\mathcal{F}}_{\mathrm{loc}}^{\ddag}$. Indeed,
the polarity of $\{0\}$ implies $C^1(\mathbb{R}^N\setminus\{0\}
)\subset\dot
{\mathcal{F}}_{\mathrm{loc}}$ and
we have that for any compact set $K(\subset\mathbb{R}^N\setminus\{0\}
)$ with
its relatively compact neighborhood $G(\subset\mathbb{R}^N\setminus\{
0\})$
satisfying $d(K,G^c):={\inf_{x\in K,y\in G^c}}|x-y|>0$
\[
\iint_{K\times G^c}\frac{(u(x)-u(y))^2}{|x-y|^{N+\alpha}}\,dx\,dy\leq
\frac{|K|\|u\|^2_{C^{0,\beta/2}}\sigma(\mathbb{S}^{N-1})}{(\alpha
-\beta)d(K,G^c)^{\alpha-\beta}}<\infty,
\]
equivalently,
\[
\iint_{K\times\mathbb{R}^N}\frac{(u(x)-u(y))^2}{|x-y|^{N+\alpha
}}\,dx\,dy<\infty,
\]
where $|K|$ is the volume of $K$ and
$\|u\|_{C^{0,\beta/2}}:={\sup_{x\ne y}}\frac{|u(x)-u(y)|}{|x-y|^{\beta/2}}$.
Therefore $u$ admits the generalized Fukushima decomposition.
\end{ex}

\section*{Acknowledgments}
This paper is based on the first part of the author's unpublished draft.
The author would like to thank Professors Zhen-Qing Chen and Masatoshi Fukushima
for their valuable comments.

%

%
\printaddresses

\end{document}


\begin{frontmatter}

\title{Errata\\
Stochastic calculus
over symmetric Markov processes without time reversal}\vspace*{6pt}
\runtitle{Errata}
\pdftitle{Errata
Stochastic calculus
over symmetric Markov processes without time reversal by K. Kuwae}

\textit{Ann. Probab.} \textbf{38} (2010)
1532--1569

\begin{aug}
\author[A]{\fnms{Kazuhiro} \snm{Kuwae}\corref{}\ead[label=e3]{kuwae@gpo.kumamoto-u.ac.jp}}
\runauthor{K. Kuwae}
\affiliation{Kumamoto University}
\address[A]{Department of Mathematics and Engineering\\
Graduate School of Science and Technology\\
Kumamoto University\\
Kumamoto, 860-8555 \\
Japan \\
\printead{e3}} 
\end{aug}

\received{\smonth{7} \syear{2011}}
\revised{\smonth{8} \syear{2011}}


%
\begin{keyword}[class=AMS]
\kwd[Primary ]{31C25}
\kwd[; secondary ]{60J25}
\kwd{60J45}
\kwd{60J75}.
\end{keyword}
\begin{keyword}
\kwd{Symmetric Markov process}
\kwd{Dirichlet form}
\kwd{Revuz measure}
\kwd{martingale additive functionals of finite energy}
\kwd{continuous additive functional of zero energy}
\kwd{Nakao's CAF of zero energy}
\kwd{Fukushima decomposition}
\kwd{semi-martingale}
\kwd{Dirichlet processes}
\kwd{stochastic integral}
\kwd{It\^{o} integral}
\kwd{Fisk--Stratonovich integral}
\kwd{time reversal operator}
\kwd{dual
predictable projection}.
\end{keyword}

\vspace*{-3pt}

\end{frontmatter}

\section{Errata}

The sentence ``In view\vspace*{1pt} of Theorem 2.2 in [22]$,\ldots$ for
q.e. \mbox{$x\in E$}.'' at page 1538 should be eliminated. The
definitions of $\mathring{\mathcal{M}}^d$, $\mathring{\mathcal{M}}^j$
and $\mathring{\mathcal{M}}^{\kappa}$ are corrected to be like
$\mathring{\mathcal{M}}^d:=\{M\in \mathring{\mathcal{M}} \mid \langle
M,N\rangle\equiv 0$ for $N\in \mathring{\mathcal{M}}^c\}$.

The statements of Theorem 2.1 and Corollary 2.3 in \cite{KwstochI} are
incorrect, which come from the error in \cite{CFKZStoch} (see
\cite{CFKZErataStoch}). The corrected statement of Theorem 2.1 in
\cite{KwstochI} can be found below. Its proof can be obtained in the
same way as in \cite{KwstochI}. The class $\widehat {\J}$ introduced in
\cite{KwstochI} is unnecessary for the corrected statement.\vspace*{-3pt}
\begin{theorem}[(Corrected statement of Theorem 2.1 in \cite
{KwstochI})]\label
{thmpurediscontimartincorrespondince}
There exists a one-to-one correspondence between ${\J}/\sim$ and
$\M_{\loc}^{d,\lbracket0,\zeta\lbracket}$
which is characterized by the relation that
for $\phi\in{\J}$ (resp., $M\in\M_{\loc}^{d,\lbracket0,\zeta\lbracket}$),
there exists $M\in\M_{\loc}^{d,\lbracket0,\zeta\lbracket}$ (resp.,
$\phi\in
{\J}$)
such that $\Delta M _t=\phi(X_{t-},X_t)$,
$t\in[0,\zeta[$, $\PP_x$-a.s. for q.e. $x\in E$.
Moreover, we have
$\langle M\rangle_t=\int_0^t\int_{E_{\partial}}\phi^2(X_s,y)N(X_s,dy)\,
dH_s$ for all
$t\in[0,\infty[$ $\PP_x$-a.s. for q.e. $x\in E$.
\end{theorem}

We define subclasses of $\M^{d,\lbracket0,\zeta\lbracket}_{\loc}$ as follows:
\begin{eqnarray*}
\M^{j,\lbracket0,\zeta\lbracket}_{\loc}&:=&\bigl\{M\in\M^{d,\lbracket
0,\zeta
\lbracket}_{\loc} \mid\phi(\cdot,\partial)=0, \kappa\mbox{-a.e. on
}E\bigr\}
, \\
\M^{\kappa,\lbracket0,\zeta\lbracket}_{\loc}&:=&\bigl\{M\in\M
^{d,\lbracket0,\zeta
\lbracket}_{\loc} \mid\phi=0, J\mbox{-a.e. on }E\times E\bigr\}.
\end{eqnarray*}
Then\vspace*{1pt} we have that
$M\in\M^{j,\lbracket0,\zeta\lbracket}_{\loc}$, $N\in\M^{\kappa
,\lbracket0,\zeta
\lbracket}_{\loc}$ imply $\langle M,N\rangle\equiv0$ $\PP_x$-a.s. for
q.e. $x\in E$, and\vspace*{1pt} every $M\in\M_{\loc}^{\lbracket0,\zeta
\lbracket}$ is
decomposed to $M=M^c+M^j+M^{\kappa}$, where
$M^c\in\M^{c,\lbracket0,\zeta\lbracket}_{\loc}$, $M^j\in
M^{j,\lbracket0,\zeta
\lbracket}_{\loc}$, $M^{\kappa}\in\M^{\kappa,\lbracket0,\zeta\lbracket
}_{\loc}$
have the properties
$\langle M^c,M^j\rangle\equiv\langle M^j,M^{\kappa}\rangle\equiv
\langle M^{\kappa},M^c\rangle\equiv
0$ in view of
Theorem \ref{thmpurediscontimartincorrespondince}.
The statement of Remark 2.3 in \cite{KwstochI} is changed to be the following:\vspace*{-3pt}
\begin{rem}
For each $i=c,d,j,\kappa$, we let
\begin{eqnarray*}
&&\accentset{\ \circ}{\M}^{i+}_{ \loc}:=\bigl\{M\mid\mbox
{ there exists }\{G_n\}\in\Theta\mbox{ and }M^{(n)}\in
\accentset{\ \circ}{\M}^i\mbox{ such that } \\
&&\hphantom{\accentset{\ \circ}{\M}^{i+}_{ \loc}:=\bigl\{}
M_t=M_t^{(n)}\mbox{ for all }t\leq\tau_{G_n}\mbox{ and }n\in
\mathbb{N}, \PP_x\mbox{-a.s. for q.e. }x\in E\bigr\}.
\end{eqnarray*}
Then $\M_{\loc}^{i,\lbracket0,\zeta\lbracket}=\accentset{\ \circ}{\M}
^{i+}_{ \loc}$ ($i=c,d,j,\kappa$).
More strongly, we have $\M_{\loc}^{c,\lbracket0,\zeta\lbracket
}=\accentset{\ \circ}{\M}^c_{ \loc}$.\vspace*{-3pt}
\end{rem}


All the statements of Corollary 4.1 and 4.3, Definition 4.3 and
Theorem~4.1 in~\cite{KwstochI}, the generalized Fukushima decomposition
(Theorem~4.2 in~\cite{KwstochI}),
the generalized It\^o formula (Theorem 4.3 in \cite{KwstochI}) and
their corollaries
(Corollaries~4.4 and~4.5 in \cite{KwstochI}) hold only for $t\in[0,
\zeta[$ $\PP_x$-a.s. for q.e. $x\in E$. The proofs of them can be
done in the same way as in \cite{KwstochI}. The classes $\dot{\F
}_{\loc}^{\ddag}$ and $\F_{\loc}^{\ddag}$ introduced in \cite
{KwstochI} are unnecessary for
the corrected statements.
We only expose the corrected statement of generalized Fukushima
decomposition below for completeness.\vspace*{-3pt}
\begin{theorem}[(Corrected statement of Theorem 4.2 in \cite
{KwstochI})]\label{thmFDlocallyinDir}
For $u\in\dot{\F}_{\loc}^{\dag}$, the additive functional $A^u$
defined by $A_t^u:=u(X_t)-u(X_0)$ can be decomposed as
\[
A^u=M^u+N^u,\qquad M^u\in\M^{\lbracket0,\zeta\lbracket}_{\loc},\qquad
N^u\in\mathcal{N}_{c,\loc}
\]
in the sense that $A_t^u=M_t^u+N_t^u, t\in[0,\zeta[$, $\PP_x\mbox
{-a.s. for q.e. }x\in E$.
Such a decomposition is unique up to the equivalence of additive
functionals on $\lbracket0,\zeta\lbracket$ (or of
local additive functionals).\vspace*{-3pt}
\end{theorem}



%

\printaddresses